\documentclass[11pt]{amsart}

\usepackage{amsmath,amssymb,ifthen,pstricks,pst-node}
\usepackage[breaklinks]{hyperref}
\usepackage[round,sort&compress]{natbib}
\usepackage[top=3cm,bottom=3cm,left=3.4cm,right=3.4cm]{geometry}

\newcounter{isamac} 
\ifthenelse{\value{isamac}=1}{
\input BoxedEPS          %
\SetTexturesEPSFSpecial  
\HideDisplacementBoxes   %
}{}

\newcommand{\an}{{\rm An}}       
\newcommand{\pa}{{\rm pa}}       
\newcommand{\un}{{\rm un}}       

\newcommand{\indm}[2]{\ensuremath{{\mathfrak I}_{\kern-1pt\scriptstyle#1}({\mathcal
#2})}} 

\newcommand{\ind}{\mbox{$\perp \kern-5.5pt \perp$}}

\newcommand{\uned}{\hbox{\kern3pt\raise2.5pt\vbox{\hrule
width9pt height 0.3pt}\kern3pt}}

\newcommand{\dashed}{\hbox{\kern3.05pt\raise2.5pt\vbox{\hrule
width1.7pt height 0.3pt}\kern1.8pt\raise2.5pt\vbox{\hrule
width1.7pt height 0.3pt}\kern1.8pt\raise2.5pt\vbox{\hrule
width1.7pt height 0.3pt}\kern1.8pt\raise2.5pt\vbox{\hrule
width1.7pt height 0.3pt}\kern3.05pt}}

\newcommand{\lhead}{\ensuremath{\prec}}

\newcommand{\head}{\ensuremath{\succ}}

\newcommand{\pedg}[2]{\ensuremath{{\kern0.5pt
\scriptstyle{\ifthenelse{\equal{\head}{#1}}{\lhead\kern0.5pt}{#1\kern0.5pt}}\joinrel\relbar
\negthinspace\relbar\joinrel{\kern0.5pt #2}\kern0.5pt}}}

\newcommand{\pdots}{\hbox{\kern2.5pt\raise1.5pt\hbox{\ensuremath{\ldots}}\ke
rn2.5pt}}  

\newcommand{\RRR}{\mathbb{R}}
\newcommand{\ND}{\mathcal{N}}
\newcommand{\tr}{\mathrm{tr}}
\newcommand{\bd}{\mathrm{bd}}
\newcommand{\Bd}{\mathrm{Bd}}

\newcommand{\bi}{\leftrightarrow}
\newcommand{\arr}{\mathrm{arr}}

\newtheorem{theorem}{Theorem}
\newtheorem{lemma}[theorem]{Lemma} 
\newtheorem{proposition}[theorem]{Proposition} 
\newtheorem{remark}[theorem]{Remark}
\newtheorem{corollary}[theorem]{Corollary}
\newtheorem{definition}[theorem]{Definition}

\title[Graphical Methods for Covariance Models]{Graphical Methods for
  Efficient Likelihood Inference in 
  Gaussian Covariance Models}

\author{Mathias Drton}
\address{Department of Statistics\\
       University of Chicago\\
       5734 S. University Ave\\
       Chicago, IL 60637, USA}
\email{drton@uchicago.edu} 
\thanks{{\em Acknowledgments.} This paper is based upon work supported by the U.S.~National Science
  Foundation (DMS-0505612, 0505865) and the U.S.~National Institutes for
  Health (R01-HG2362-3).}
\author{Thomas S.\ Richardson}
\address{Department of Statistics\\
       University of Washington\\
       Seattle, WA 98195-4322, USA}
\email{tsr@stat.washington.edu}

\begin{document}

\begin{abstract}
  In graphical modelling, a bi-directed graph encodes marginal
  independences among random variables that are identified with the
  vertices of the graph.  We show how to transform a bi-directed graph into
  a maximal ancestral graph that (i) represents the same independence
  structure as the original bi-directed graph, and (ii) minimizes the
  number of arrowheads among all ancestral graphs satisfying (i). Here the
  number of arrowheads of an ancestral graph is the number of directed
  edges plus twice the number of bi-directed edges. In Gaussian models,
  this construction can be used for more efficient iterative maximization
  of the likelihood function and to determine when maximum likelihood
  estimates are equal to empirical counterparts.
\end{abstract}


\maketitle

\section{Introduction}\label{sec:introduction}

\begin{figure}[tbp]
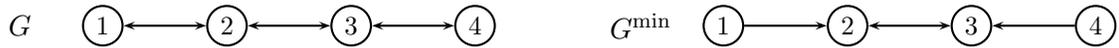

  \centering
  \begin{tabular}{l@{\hspace{2.25cm}}l}
    \parbox{6cm}{
      \vspace{1.1cm}   
      \small\psset{unit=1.1cm}
      \newcommand{\myNode}[2]{\circlenode{#1}{\makebox{#2}}}
      \rput(0,0.5){\normalsize $G$}
      \rput(1,0.5){\myNode{x}{$1$}} 
      \rput(2.5,0.5){\myNode{v}{$2$}} 
      \rput(4,0.5){\myNode{w}{$3$}} 
      \rput(5.5,0.5){\myNode{y}{$4$}} 
      \ncline{<->}{v}{w}
      \ncline{<->}{w}{y}
      \ncline{<->}{v}{x}
      \vspace{0.05cm}
    }
    &
    \parbox{6cm}{
      \vspace{1.1cm}   
      \small\psset{unit=1.1cm}
      \newcommand{\myNode}[2]{\circlenode{#1}{\makebox{#2}}}
      \rput(0,0.5){\normalsize $G^{\min}$}
      \rput(1,0.5){\myNode{x}{$1$}} 
      \rput(2.5,0.5){\myNode{v}{$2$}} 
      \rput(4,0.5){\myNode{w}{$3$}} 
      \rput(5.5,0.5){\myNode{y}{$4$}} 
      \ncline{<->}{v}{w}
      \ncline{<-}{w}{y}
      \ncline{<-}{v}{x}
      \vspace{0.05cm}
    }
  \end{tabular}
  \caption{A bi-directed graph $G$ with (unique) minimally oriented graph $G^{\min}$.}
  \label{fig:intro}
\end{figure}

In graphical modelling, bi-directed graphs encode marginal independences
among random variables that are identified with the vertices of the graph
\citep{pearl:1994,kauermann:dual, richardson:2003}.  In particular, if two
vertices are not joined by an edge, then the two associated random
variables are assumed to be marginally independent.  For example, the graph
$G$ in Figure \ref{fig:intro}, whose vertices are to be identified with a
random vector $(X_1,X_2,X_3,X_4)$, represents the pairwise marginal
independences $X_1\ind X_3$, $X_1\ind X_4$, and $X_2\ind X_4$.  While other
authors \citep{coxwerm:lindep,coxwerm:book,edwards:2000} have used dashed
edges to represent marginal independences, the bi-directed graphs we employ
here make explicit the connection to path diagrams
\citep{wright:methodpathcoeff,koster:pathdiagram}.

Gaussian graphical models for marginal independence, also known as
covariance graph models, impose zero patterns in the covariance matrix,
which are linear hypotheses on the covariance matrix \citep{anderson:1973}.
The graph in Figure \ref{fig:intro}, for example, imposes
$\sigma_{13}=\sigma_{14}=\sigma_{24}=0$.  An estimation procedure designed
for covariance graph models is described in \cite{drton:2003b}; see also
\cite{chaudhuri:2007}.  Other recent work involving these models includes
\cite{mao:2004} and \cite{wermuth:2006}.

In this paper we employ the connection between bi-directed graphs and the
more general ancestral graphs with undirected, directed, and bi-directed
edges (Section \ref{sec:graphterm}).  For the statistical motivation of
ancestral graphs see \cite{richardson:2002}; for causal interpretation see
\cite{richardsonspirtes:2003}. We show how to construct a maximal ancestral
graph $G^{\min}$, which we call a minimally oriented graph, that is Markov
equivalent to a given bi-directed graph $G$ and such that the number of
arrowheads is minimal (Sections
\ref{sec:simpl-sets}--\ref{sec:orient-edges}).  Two ancestral graphs are
Markov equivalent if the independence models associated with the two graphs
coincide; see for example \cite{roverato:2005} for some recent results on Markov
equivalence of different types of graphs.  The number of arrowheads is the
number of directed edges plus twice the number of bi-directed edges.
Minimally oriented graphs provide useful nonparametric information about
Markov equivalence of bi-directed, undirected and directed acyclic graphs.
For example, the graph $G$ in Figure \ref{fig:intro} is not Markov
equivalent to an undirected graph because $G^{\min}$ is not an undirected
graph, and $G$ is not Markov equivalent to a DAG because $G^{\min}$
contains a bi-directed edge.  The graph in Figure \ref{fig:intro} has a
unique minimally oriented graph but in general, minimally oriented graphs
are not unique.  Our construction procedure
(Algorithm~\ref{alg:construct-Gos}) involves a choice of a total order
among the vertices.  Varying the order one may obtain all minimally
oriented graphs.

For covariance graph models, minimally oriented graphs allow one to
determine when the maximum likelihood estimate of a variance or covariance
is available explicitly as its empirical counterpart (Section
\ref{sec:impl-max-lik}).  For example, since no arrowheads appear at the
vertices $1$ and $4$ in the graph $G^{\min}$ in Figure \ref{fig:intro}, the
maximum likelihood estimates of $\sigma_{11}$ and $\sigma_{44}$ must be
equal to the empirical variance of $X_1$ and $X_4$, respectively.  The
likelihood function for covariance graph models may be multi-modal, though
simulations suggest this only occurs at small sample sizes, or under
mis-specification \citep{drton:2004}.  However, when a minimally oriented
graph reveals that a parameter estimate is equal to an empirical quantity
(such as $\sigma_{11}$ and $\sigma_{44}$ in the above example) then even if
the likelihood function is multi-modal this parameter will take the same
value at every mode.  Perhaps most importantly, minimally oriented graphs
allow for computationally more efficient maximum likelihood fitting; see
Remark \ref{rem:efficient} and the example in Section \ref{sec:example}.

\section{Ancestral graphs and their global Markov property}\label{sec:graphterm} 

This paper deals with {\em simple mixed graphs\/}, which feature undirected
($v-w$), directed ($v\to w$) and bi-directed edges ($v\bi w$) under the
constraint that there is at most one edge between two vertices.  In this
section we give a formal definition of these graphs and discuss their
Markov interpretation.

\subsection{Simple mixed graphs}
\label{subsec:simple-mixed-graphs}

Let $\mathcal{E}=\{\emptyset, - ,\leftarrow,\to,\bi\}$ be the set of
possible edges between an ordered pair of vertices; $\emptyset$ denoting
that there is no edge.  A {\em simple mixed graph\/} $G=(V,E)$ is a pair of
a finite {\em vertex set\/} $V$ and an {\em edge map\/} $E:V\times V\to
\mathcal{E}$.  The edge map $E$ has to satisfy that for all $v,w\in V$,
\begin{enumerate}
\item $E(v,v)=\emptyset$, i.e., there is no edge between a vertex and
  itself,
\item $E(v,w)=E(w,v)$ if $E(v,w)\in\{-,\bi\}$,
\item $E(v,w)=\,\to\;\; \iff\;\;  E(w,v)=\,\leftarrow$.
\end{enumerate}
In the sequel, we write $v-w\in G$, $v\to w\in G$, $v\leftarrow w\in G$ or
$v\bi w\in G$ if $E(v,w)$ equals $-$, $\to$, $\leftarrow$ or $\bi$,
respectively.  If $E(v,w)\not= \emptyset$, then $v$ and $w$ are {\em
  adjacent\/}.  If there is an edge $v\leftarrow w\in G$ or $v\bi w\in G$
then there is an {\em arrowhead at $v$\/} on this edge. If there is an edge
$v\to w\in G$ or $v - w\in G$ then there is a {\em tail at $v$\/} on this
edge.  A vertex $w$ is in the {\em boundary\/} of $v$, denoted by $\bd(v)$,
if $v$ and $w$ are adjacent.  The boundary of vertex set $A\subseteq V$ is
the set $\bd(A) = [\cup_{v\in A} \bd(v)] \setminus A$.  We write
$\Bd(v)=\bd(v)\cup \{v\}$ and $\Bd(A)=\bd(A)\cup A$.  An induced subgraph
of $G$ over a vertex set $A$ is the mixed graph $G_A=(A,E_A)$ where $E_A$
is the restriction of the edge map $E$ on $A\times A$.  The {\em
  skeleton\/} of a simple mixed graph is obtained by making all edges
undirected.

In a simple mixed graph a sequence of adjacent vertices $(v_1,\ldots,v_k)$
uniquely determines the sequence of edges joining consecutive vertices
$v_i$ and $v_{i+1}$, $1\le i\le k-1$.  Hence, we can define a {\em path\/}
$\pi$ between two vertices $v$ and $w$ as a sequence of distinct vertices
$\pi=(v,v_1,\ldots,v_k,w)$ such that each vertex in the sequence is
adjacent to its predecessor and its successor.  A path $v\to\cdots\to w$
with all edges of the form $\to$ and pointing toward $w$ is a directed path
from $v$ to $w$.  If there is such a directed path from $v$ to $w\not= v$,
or if $v=w$, then $v$ is an {\em ancestor\/} of $w$.  We denote the set of
all ancestors of a vertex $v$ by $\an(v)$ and for a vertex set $A\subseteq
V$ we define $\an(A)=\cup_{v\in A} \an(v)$.  Finally, a directed path from
$v$ to $w$ together with an edge $w\to v\in G$ is called a {\em directed
  cycle\/}.

Important subclasses of simple mixed graphs are illustrated in Figure
\ref{fig:exsimplemixed}.  {\em Bi-directed\/}, {\em undirected\/} and {\em
  directed graphs\/} contain only one type of edge.  {\em Directed acyclic
  graphs\/} (DAGs) are directed graphs without directed cycles. These three
types of graphs are special cases of {\em ancestral graphs\/}
\citep{richardson:2002}.
\begin{definition}
  \label{def:ancestralgraph} A simple mixed graph $G$ is an {\em ancestral
    graph\/} if it holds that
\begin{enumerate}
\item $G$ does not contain any directed cycles;
\item if $v-w\in G$, then there does not exist $u$ such that $u\to v\in G$
  or $u\bi v\in G$;
\item if $v\bi w\in G$, then $v$ is not an ancestor of $w$.
\end{enumerate}
\end{definition}

\subsection{Global Markov property for ancestral graphs}
 \label{subsec:global-markov}

Ancestral graphs can be given an independence interpretation, known as the
global Markov property, by a graphical separation criterion called
$m$-separation \cite[\S3.4]{richardson:2002}.  An extension of Pearl's
(\citeyear{pearl:1988}) $d$-separation for DAGs, $m$-separation uses the
notion of {\em colliders\/}.  A non-endpoint vertex $v_i$ on a path is a
{\em collider on the path\/} if the edges preceding and succeeding $v_i$ on
the path both have an arrowhead at $v_i$, that is, $v_{i-1}\to
v_i\leftarrow v_{i+1}$, $v_{i-1}\to v_i\bi v_{i+1}$, $v_{i-1}\bi
v_i\leftarrow v_{i+1}$ or $v_{i-1}\bi v_i\bi v_{i+1}$ is part of the path.
A non-endpoint vertex that is not a collider is a {\em non-collider on the
  path\/}.
\begin{definition} 
  \label{def:m-connect} 
  A path $\pi$ between vertices $v$ and $w$ in a simple mixed graph $G$ is
  {\em $m$-connecting\/} given a possibly empty set $C\subseteq
  V\setminus\{v,w\}$ if (i) every non-collider on $\pi$ is not in $C$, and
  (ii) every collider on $\pi$ is in $\an(C)$.  If no path $m$-connects $v$
  and $w$ given $C$, then $v$ and $w$ are {\em $m$-separated\/} given $C$.
  Two non-empty and disjoint sets $A$ and $B$ are $m$-separated given
  $C\subseteq V\setminus (A\cup B)$, if any two vertices $v\in A$ and $w\in
  B$ are $m$-separated given $C$.
\end{definition}

\begin{figure}[tbp]
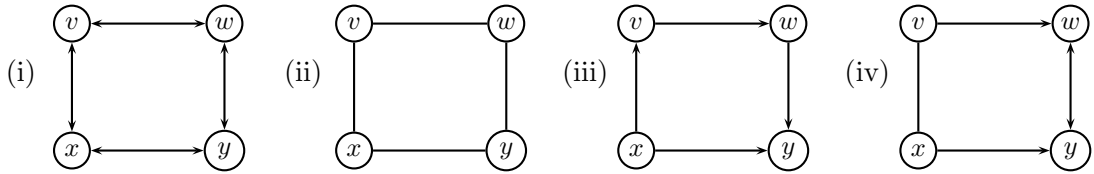

  \centering
  \begin{tabular}{l@{\hspace{0cm}}l@{\hspace{0cm}}l@{\hspace{0cm}}l}
    \hspace{-0.75cm}
    \parbox{3.75cm}{
      \vspace{2.75cm}   
      \small\psset{unit=1.35cm}
      \newlength{\MyLength}
      \settowidth{\MyLength}{$v$}
      \newcommand{\myNode}[2]{\circlenode{#1}{\makebox[\MyLength]{#2}}}
      \rput(0.5,1.25){\normalsize (i)}
      \rput(1,1.75){\myNode{v}{$v$}} 
      \rput(1,0.5){\myNode{x}{$x$}} 
      \rput(2.5,1.75){\myNode{w}{$w$}}
      \rput(2.5,0.5){\myNode{y}{$y$}}
      \ncline{<->}{v}{w}
      \ncline{<->}{x}{y}
      \ncline{<->}{w}{y}
      \ncline{<->}{v}{x}
    }
    &
    \parbox{3.75cm}{
      \vspace{2.75cm}   
      \small\psset{unit=1.35cm}
      \settowidth{\MyLength}{$v$}
      \newcommand{\myNode}[2]{\circlenode{#1}{\makebox[\MyLength]{#2}}}
      \rput(0.5,1.25){\normalsize (ii)}
      \rput(1,1.75){\myNode{v}{$v$}} 
      \rput(1,0.5){\myNode{x}{$x$}} 
      \rput(2.5,1.75){\myNode{w}{$w$}}
      \rput(2.5,0.5){\myNode{y}{$y$}}
      \ncline{-}{v}{w}
      \ncline{-}{x}{y}
      \ncline{-}{w}{y}
      \ncline{-}{v}{x}
    }    
    &
    \parbox{3.75cm}{
      \vspace{2.75cm}   
      \small\psset{unit=1.35cm}
      \settowidth{\MyLength}{$v$}
      \newcommand{\myNode}[2]{\circlenode{#1}{\makebox[\MyLength]{#2}}}
      \rput(0.5,1.25){\normalsize (iii)}
      \rput(1,1.75){\myNode{v}{$v$}} 
      \rput(1,0.5){\myNode{x}{$x$}} 
      \rput(2.5,1.75){\myNode{w}{$w$}}
      \rput(2.5,0.5){\myNode{y}{$y$}}
      \ncline{->}{v}{w}
      \ncline{->}{x}{y}
      \ncline{->}{w}{y}
      \ncline{<-}{v}{x}
    }    
    &
    \parbox{3.75cm}{
      \vspace{2.75cm}   
      \small\psset{unit=1.35cm}
      \settowidth{\MyLength}{$v$}
      \newcommand{\myNode}[2]{\circlenode{#1}{\makebox[\MyLength]{#2}}}
      \rput(0.5,1.25){\normalsize (iv)}
      \rput(1,1.75){\myNode{v}{$v$}} 
      \rput(1,0.5){\myNode{x}{$x$}} 
      \rput(2.5,1.75){\myNode{w}{$w$}}
      \rput(2.5,0.5){\myNode{y}{$y$}}
      \ncline{->}{v}{w}
      \ncline{->}{x}{y}
      \ncline{<->}{w}{y}
      \ncline{-}{v}{x}
    }
  \end{tabular}
  \caption{Simple mixed graphs. (i) A bi-directed graph,
    (ii) an undirected graph, (iii) a DAG, (iv) an ancestral graph.}
  \label{fig:exsimplemixed}
\end{figure}

Let $G=(V,E)$ be an ancestral graph whose vertices index a random vector
$(X_v\mid v\in V)$.
For $A\subseteq V$, let $X_A$ be the subvector $(X_v\mid v\in A)$. The {\em
  global Markov property\/} for $G$ states that $X_A$ is conditionally
independent of $X_B$ given $X_C$ whenever $A$, $B$ and $C$ are pairwise
disjoint subsets such that $A$ and $B$ are $m$-separated given $C$ in $G$.
Subsequently, we denote such conditional independence using the shorthand
$A\ind B\mid C$ that avoids making the probabilistic context explicit.  The
global Markov property, when applied to each of the graphs in Figure
\ref{fig:exsimplemixed} in turn, implies (among other independences) that:
\begin{enumerate}
\item  $v\ind y$ and $w\ind x$;
\item  $v\ind y\mid \{w,x\}$ and
$w\ind x\mid \{v,y\}$;
\item  $v\ind y\mid \{w,x\}$ and $w\ind x\mid v$;
\item $v\ind y\mid x$ and $w\ind x\mid v$.
\end{enumerate}

If $G$ is a bi-directed graph, then the global Markov property states the
marginal independence $v\ind w$ if $v$ and $w$ are not adjacent.  In a
multivariate normal distribution such pairwise marginal independences hold
iff all independences stated by the global Markov property for $G$ hold
\citep{kauermann:dual}.  Without any distributional assumption,
\citet[\S4]{richardson:2003} shows that the independences stated by the
global Markov property of a bi-directed graph hold iff certain (not only
pairwise) marginal independences hold; see also \cite{matus:1994}.

The graphs in Figure \ref{fig:exsimplemixed} have the property that for
every pair of non-adjacent vertices $v$ and $w$ there exists some subset
$C$ such that the global Markov property states that $v\ind w\mid C$.
Ancestral graphs with this property are called {\em maximal\/}.  If an
ancestral graph $G$ is not maximal, then there exists a unique Markov
equivalent maximal ancestral graph $\bar G$ that contains all the edges
present in $G$.  Moreover, any edge in $\bar G$ that is not present in $G$
is bi-directed \cite[][\S3.7]{richardson:2002}.  Two ancestral graphs $G_1$
and $G_2$ are {\em Markov equivalent\/} if they have the same vertex set
and the global Markov property states the same independences for $G_1$ as
for $G_2$.

The following facts are easily established; see also \cite{richardson:2002}.

\begin{lemma} 
  \label{lem:MAGs} 
  \begin{enumerate}
  \item  Markov equivalent maximal ancestral graphs have the same skeleton.
  \item If $\bar G$ is an ancestral graph that is Markov equivalent to a
    maximal ancestral graph $G$ and has the same skeleton as $G$, then
    $\bar G$ is also a maximal ancestral graph.
  \item Bi-directed, undirected and directed acyclic graphs are maximal
    ancestral graphs.
  \end{enumerate}
\end{lemma}

\subsection{Boundary containment}
\label{subsec:boundary-contain}

In the subsequent Sections~\ref{sec:simpl-sets} and \ref{sec:orient-edges}
we will construct maximal ancestral graphs that are Markov equivalent to a
given bi-directed graph.  Via Theorem~\ref{thm:get-mconnect} below, the
following property plays a crucial role in these constructions.

\begin{definition}
  \label{def:dir-bound-contain}
  A simple mixed graph $G$ has the {\em boundary containment property} if
  for all distinct vertices $v,w\in V$ the presence of an edge $v-w$
  implies that $\Bd(v)=\Bd(w)$ and the presence of an edge $v\to w$ in $G$
  implies that $\Bd(v)\subseteq \Bd(w)$.
\end{definition}

In the Appendix we present lemmas on the structure of $m$-connecting paths
in graphs with the boundary containment property.  These lemmas yield the
following key result.

\begin{theorem}
  \label{thm:get-mconnect}
  If $\bar G$ is an ancestral graph that has the same skeleton as a bi-directed
  graph $G$, then $G$ and $\bar G$ are Markov equivalent iff $\bar G$ has
  the boundary containment property.
\end{theorem}
\begin{proof}
  Two vertices are adjacent in $G$ iff they are adjacent in $\bar G$.
  Therefore, $G$ and $\bar G$ are Markov equivalent iff it holds that two
  non-adjacent vertices $v$ and $w$ are $m$-connected given $C\subseteq V$
  in $G$ iff they are $m$-connected given $C$ in $\bar G$.
  
  ($\Longrightarrow:$) Suppose $\bar G$ does not have the boundary
  containment property, i.e., there exists an edge $v-w\in \bar G$ or an
  edge $v\to w\in\bar G$ such that $\Bd(v)\not\subseteq\Bd(w)$.  Choose
  $u\in\Bd(v)\setminus\Bd(w)$.  Since $u$ and $w$ are not adjacent, they
  are $m$-separated given $C=\emptyset$ in $G$.  In $\bar G$, however, the
  path $(u,v,w)$ $m$-connects $u$ and $w$ given $C=\emptyset$.  Hence, $G$
  and $\bar G$ are not Markov equivalent.
  
  ($\Longleftarrow:$) First, let $v$ and $w$ be non-adjacent vertices that
  are $m$-connected given $C\subseteq V$ in $\bar G$.  By
  Lemma~\ref{prop:connect6}, there is a path $\bar\pi = (v,v_1,\ldots,v_k,w)$
  that $m$-connects $v$ and $w$ given $C$ in $\bar G$ and is such that
  $v_1,\ldots,v_k$ are colliders with $\{v_1\ldots,v_k\}\subseteq C$.
  Since $G$ is a bi-directed graph, the corresponding path $\pi=
  (v,v_1,\ldots,v_k,w)$ in $G$ also $m$-connects $v$ and $w$ given $C$.
  
  Conversely, let $v$ and $w$ be non-adjacent vertices that are
  $m$-connected given $C\subseteq V$ in $G$.  Let $\pi =
  (v_0,v_1,\ldots,v_k,v_{k+1})$ $m$-connect $v=v_0$ and $w=v_{k+1}$ given
  $C$ in $G$ such that no shorter path $m$-connects $v$ and $w$ given $C$.
  Then $v_1,\ldots,v_k$ are colliders, $\{v_1\ldots,v_k\}\subseteq C$, and
  $v_{i-1}$ and $v_{i+1}$, $i=1,\ldots,k$, are not adjacent in $G$.  (This
  is a special case of Lemmas~\ref{lem:tsr2} and \ref{prop:connect6}
  because a bi-directed graph trivially satisfies the boundary containment
  property.)  It follows that, for all $i=1,\ldots,k-1$,
  $v_{i-1}\in\Bd(v_i)$ but $v_{i-1}\notin\Bd(v_{i+1})$, and similarly
  $v_{i+2}\notin\Bd(v_i)$ but $v_{i+2}\in\Bd(v_{i+1})$.  This implies that
  $\Bd(v_i)\not\subseteq \Bd(v_{i+1})$ and $\Bd(v_i)\not\supseteq
  \Bd(v_{i+1})$ for all $i=1,\ldots,k-1$.  Since $\bar G$ has the boundary
  containment property, it must hold that $v_i\bi v_{i+1}\in \bar G$ for
  all $i=1,\ldots,k-1$.  Therefore, $v_2,\ldots,v_{k-1}$ are colliders on
  the path $\bar\pi=(v,v_1,\ldots,v_k,w)$ in $\bar G$.  Similarly, it
  follows that $v_2\in\Bd(v_1)\setminus\Bd(v)$, which entails
  $\Bd(v_1)\not\subseteq \Bd(v)$.  Thus, $v_1$ is a collider on $\bar\pi$.
  Analogously, we can show that $v_k$ is a collider on $\bar\pi$, which
  yields that $\bar\pi$ is a path in $\bar G$ that $m$-connects $v$ and $w$
  given $C$.
\end{proof}

\section{Simplicial graphs}\label{sec:simpl-sets}

In this section we show how simplicial vertex sets of a bi-directed graph
can be used to construct a Markov equivalent maximal ancestral graph by
removing arrowheads from certain bi-directed edges.  Simplicial sets are
also important in other contexts such as collapsibility \citep[\S2.1.3,
p.121 and 219]{madigan:1990,kauermann:dual,lau:bk} and triangulation of
graphs \citep[\S5.3]{jensen:2001}.

\begin{definition} 
  \label{def:simplicialset}
  A vertex $v\in V$ is {\em simplicial\/}, if $\Bd(v)$ is
  complete, i.e., every pair of vertices in $\Bd(v)$ are adjacent. Similarly, a set 
  $A\subseteq V$ is {\em simplicial\/}, if $\Bd(A)$ is complete.  
\end{definition}

Simplicial vertices can be characterized in terms of boundary containment
as follows.

\begin{proposition}\label{prop:completebd}
  A vertex $v\in V$ is simplicial iff $\Bd (v) \subseteq \Bd (w)$ for all
  $w\in\Bd (v)$. 
\end{proposition}

If an edge between $v$ and $w$ has an arrowhead at $v$, then we say that we
{\em drop the arrowhead at $v$\/} when either $v\leftarrow w$ is replaced
by $v-w$ or $v\bi w$ is replaced by $v\to w$.

\begin{figure}[tbp]
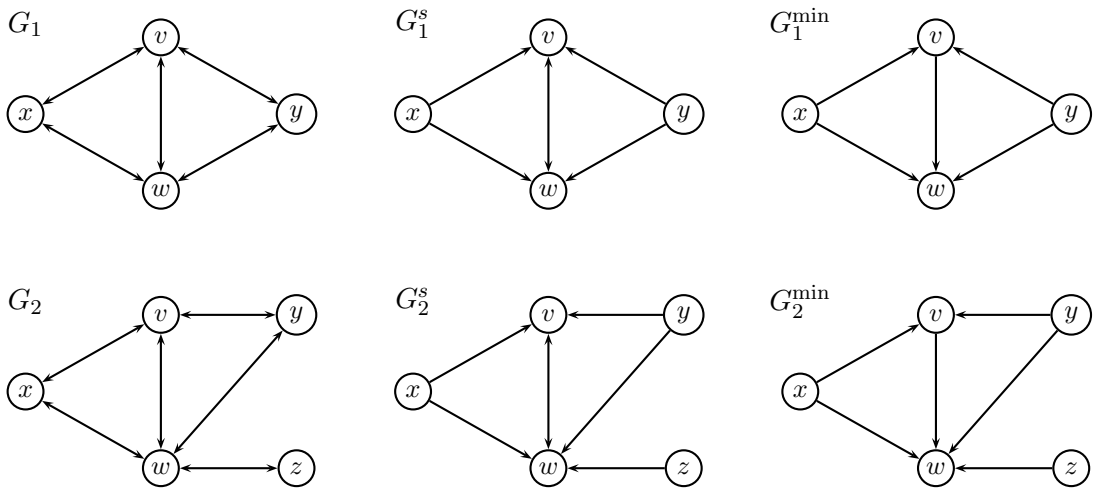

  \centering
  \hspace{-2.25cm}
  \begin{tabular}{l@{\hspace{1.15cm}}l@{\hspace{1.15cm}}l}
    \parbox{4cm}{
      \vspace{3.35cm}   
      \small\psset{unit=1.2cm}
      \settowidth{\MyLength}{$v$}
      \newcommand{\myNode}[2]{\circlenode{#1}{\makebox[\MyLength]{#2}}}
      \rput(1,2.5){\normalsize $G_1$}
      \rput(2.5,2.35){\myNode{v}{$v$}} 
      \rput(1,1.5){\myNode{x}{$x$}} 
      \rput(2.5,0.65){\myNode{w}{$w$}}
      \rput(4,1.5){\myNode{y}{$y$}}
      \ncline{<->}{v}{w}
      \ncline{<->}{w}{y}
      \ncline{<->}{v}{y}
      \ncline{<->}{w}{x}
      \ncline{<->}{v}{x}
      \vspace{0.05cm}
    }
    & 
    \parbox{4cm}{
      \vspace{3.35cm}   
      \small\psset{unit=1.2cm}
      \settowidth{\MyLength}{$v$}
      \newcommand{\myNode}[2]{\circlenode{#1}{\makebox[\MyLength]{#2}}}
      \rput(1,2.5){\normalsize $G_1^s$}
      \rput(2.5,2.35){\myNode{v}{$v$}} 
      \rput(1,1.5){\myNode{x}{$x$}} 
      \rput(2.5,0.65){\myNode{w}{$w$}}
      \rput(4,1.5){\myNode{y}{$y$}}
      \ncline{<->}{v}{w}
      \ncline{<-}{w}{y}
      \ncline{<-}{v}{y}
      \ncline{<-}{w}{x}
      \ncline{<-}{v}{x}
      \vspace{0.05cm}
    }
    & 
    \parbox{4cm}{
      \vspace{3.35cm}   
      \small\psset{unit=1.2cm}
      \settowidth{\MyLength}{$v$}
      \newcommand{\myNode}[2]{\circlenode{#1}{\makebox[\MyLength]{#2}}}
      \rput(1,2.5){\normalsize $G_1^{\min}$}
      \rput(2.5,2.35){\myNode{v}{$v$}} 
      \rput(1,1.5){\myNode{x}{$x$}} 
      \rput(2.5,0.65){\myNode{w}{$w$}}
      \rput(4,1.5){\myNode{y}{$y$}}
      \ncline{->}{v}{w}
      \ncline{<-}{w}{y}
      \ncline{<-}{v}{y}
      \ncline{<-}{w}{x}
      \ncline{<-}{v}{x}
      \vspace{0.05cm}
    }
  \\[1.75cm]
    \parbox{4cm}{
      \vspace{3.35cm}   
      \small\psset{unit=1.2cm}
      \settowidth{\MyLength}{$v$}
      \newcommand{\myNode}[2]{\circlenode{#1}{\makebox[\MyLength]{#2}}}
      \rput(1,2.5){\normalsize $G_2$}
      \rput(2.5,2.35){\myNode{v}{$v$}} 
      \rput(1,1.5){\myNode{x}{$x$}} 
      \rput(2.5,0.65){\myNode{w}{$w$}}
      \rput(4,2.35){\myNode{y}{$y$}}
      \rput(4,0.65){\myNode{z}{$z$}}
      \ncline{<->}{v}{w}
      \ncline{<->}{w}{z}
      \ncline{<->}{w}{y}
      \ncline{<->}{v}{y}
      \ncline{<->}{w}{x}
      \ncline{<->}{v}{x}
      \vspace{-0.1cm}
    }
    & 
    \parbox{4cm}{
      \vspace{3.35cm}   
      \small\psset{unit=1.2cm}
      \settowidth{\MyLength}{$v$}
      \newcommand{\myNode}[2]{\circlenode{#1}{\makebox[\MyLength]{#2}}}
      \rput(1,2.5){\normalsize $G_2^s$}
      \rput(2.5,2.35){\myNode{v}{$v$}} 
      \rput(1,1.5){\myNode{x}{$x$}} 
      \rput(2.5,0.65){\myNode{w}{$w$}}
      \rput(4,2.35){\myNode{y}{$y$}}
      \rput(4,0.65){\myNode{z}{$z$}}
      \ncline{<->}{v}{w}
      \ncline{<-}{w}{z}
      \ncline{<-}{w}{y}
      \ncline{<-}{v}{y}
      \ncline{<-}{w}{x}
      \ncline{<-}{v}{x}
      \vspace{-0.1cm}
    }
    & 
    \parbox{4cm}{
      \vspace{3.35cm}   
      \small\psset{unit=1.2cm}
      \settowidth{\MyLength}{$v$}
      \newcommand{\myNode}[2]{\circlenode{#1}{\makebox[\MyLength]{#2}}}
      \rput(1,2.5){\normalsize $G_2^{\min}$}
      \rput(2.5,2.35){\myNode{v}{$v$}} 
      \rput(1,1.5){\myNode{x}{$x$}} 
      \rput(2.5,0.65){\myNode{w}{$w$}}
      \rput(4,2.35){\myNode{y}{$y$}}
      \rput(4,0.65){\myNode{z}{$z$}}
      \ncline{->}{v}{w}
      \ncline{<-}{w}{z}
      \ncline{<-}{w}{y}
      \ncline{<-}{v}{y}
      \ncline{<-}{w}{x}
      \ncline{<-}{v}{x}
      \vspace{-0.1cm}
    }
  \end{tabular}
  \caption{Bi-directed graphs with simplicial and minimally oriented
    graphs.}
  \label{fig:Gs}
\end{figure}

\begin{definition} 
  Let $G$ be a bi-directed graph. The {\em simplicial graph\/} $G^s$ is the
  simple mixed graph obtained by dropping all the arrowheads at simplicial
  vertices of $G$.
\end{definition}
For the graph from Figure \ref{fig:intro}, $G^s$ is equal to the depicted
graph $G^{\min}$; additional examples are given in Figure \ref{fig:Gs}.
Parts (i) and (ii) of the next lemma show that simplicial graphs have
the boundary containment property. 

\begin{lemma}  
  \label{lem:orientedgesinGs}
  Let $v$ and $w$ be adjacent vertices in a simplicial graph $G^s$.  Then
  \begin{enumerate}
  \item if $v-w\in G^s$, then $\Bd(v)=\Bd(w)$;
  \item if $v\to w\in G^s$, then $\Bd(v)\subsetneq\Bd(w)$; 
  \item if $v\bi w\in G^s$, then each of $\Bd(v)=\Bd(w)$,
    $\Bd(v)\subsetneq\Bd(w)$, and
    $\Bd(v)\not\subseteq\Bd(w)\not\subseteq\Bd(v)$ might be the case.
  \end{enumerate}
\end{lemma}
\begin{proof}
  (i) and (ii) follow from Proposition \ref{prop:completebd}.  For (iii)
  see, respectively, the graphs $G_1^s$, $G_2^s$ in Figure \ref{fig:Gs},
  and $G^s=G^{\min}$ in Figure \ref{fig:intro}.
\end{proof}

\begin{theorem}  
\label{thm:equivsimple}
The simplicial graph $G^s$ of a bi-directed graph $G$ is a maximal
ancestral graph that is Markov equivalent to $G$.
\end{theorem}
\begin{proof} 
  By Lemma~\ref{lem:MAGs}, Theorem~\ref{thm:get-mconnect} and
  Lemma~\ref{lem:orientedgesinGs}, it suffices to show that $G^s$ is an
  ancestral graph.  This, however, follows from
  Lemma~\ref{lem:drop-at-simplicial} below.
\end{proof}

\begin{lemma}
  \label{lem:drop-at-simplicial}
  If $G$ is an ancestral graph that has the boundary containment property,
  then dropping all arrowheads at simplicial vertices of $G$ yields an
  ancestral graph.
\end{lemma}
\begin{proof}
  Let $\bar G$ be the graph obtained by dropping the arrowheads at
  simplicial vertices.  First, suppose $v\to w\in\bar G$ or $v\bi w\in\bar
  G$ but that there is a path $\pi$ from $w$ to $v$ that is a directed path
  in $\bar G$.  Since there are no arrowheads at simplicial vertices in
  $\bar G$, no vertex on $\pi$ including the endpoints $v$ and $w$ can be
  simplicial.  This implies that $\pi$ is a directed path from $w$ to $v$
  in $G$.  However, since $v\to w\in G$ or $v\bi w\in G$, this is a
  contradiction to $G$ being ancestral.  We conclude that $\bar G$
  satisfies conditions (i) and (iii) of
  Definition~\ref{def:ancestralgraph}.
  
  Next, suppose $v-w\in\bar G$ but that there exists another vertex $u$
  such that $u\to v\in\bar G$ or $u\bi v\in\bar G$.  It follows that $v$ is
  not simplicial.  Since $G$ is ancestral, this implies that $v\to w\in G$
  which in turn implies that $\Bd(v)\subseteq\Bd(w)$ because $G$ has the
  boundary containment property.  The set $\Bd(v)$ is not complete because
  $v$ is not simplicial.  Thus $\Bd(w)$ is not complete, i.e., $w$ is not a
  simplicial vertex.  However, this is a contradiction to the fact that
  $v\to w\in G$ but $v-w\in\bar G$.  Thus, $\bar G$ is indeed an ancestral
  graph.
\end{proof}

\begin{proposition}
\label{prop:equivsimple}
A bi-directed graph $G$ is Markov equivalent to an undirected graph iff the
simplicial graph $G^s$ induced by $G$ is an undirected graph iff $G$ is a
disjoint union of complete (bi-directed) graphs.
\end{proposition}
\begin{proof}
  If $G^s$ is an undirected graph, then by Theorem \ref{thm:equivsimple},
  $G$ is Markov equivalent to an undirected graph, namely $G^s$.
  Conversely, assume that there exists an undirected graph $U$ that is
  Markov equivalent to $G$.  Necessarily, $G$ and $U$ have the same
  skeleton (recall Lemma~\ref{lem:MAGs}).  By
  Theorem~\ref{thm:get-mconnect}, $U$ has the boundary containment
  property, which implies that every vertex is simplicial and thus that
  $G^s$ is an undirected graph (and equal to $U$).
  
  The simplicial graph $G^s$ is an undirected graph iff the vertex set of
  the inducing bi-directed graph $G$ can be partitioned into pairwise
  disjoint sets $A_1,\ldots, A_q$ such that (a) if $v\in A_i$, $1\le i\le
  q$, and $w\in A_j$, $1\le j\le q$, are adjacent, then $i=j$, and (b) all
  the induced subgraphs $G_{A_i}$, $i=1,\ldots,q$ are complete graphs
  \citep{kauermann:dual}.
\end{proof}

Under multivariate normality, a bi-directed graph that is Markov equivalent
to an undirected graph represents a hypothesis that is linear in the
covariance matrix as well as in its inverse.  The general structure of such
models is studied in \cite{jensen:1988}.

\section{Minimally oriented graphs}
\label{sec:orient-edges}

The simplicial graph $G^s$ sometimes may be a DAG.  For example, the graph
$u\bi v\bi w$ has the simplicial graph $u\to v\leftarrow w$.  However,
there exist bi-directed graphs that are Markov equivalent to a DAG and yet
the simplicial graph contains bi-directed edges.  For example, the graph
$G_1$ in Figure \ref{fig:Gs} is Markov equivalent to the DAG $G_1^{\min}$
in the same Figure.  Hence, some arrowheads may be dropped from bi-directed
edges in a simplicial graph while preserving Markov equivalence.  In this
section we construct maximal ancestral graphs from which no arrowheads may
be dropped without destroying Markov equivalence.

\subsection{Definition and construction}
\label{subsec:mini-orient-def-cons}

The following definition introduces the key object of this section.

\begin{definition}  
  \label{def:Gos}
  Let $G$ be a bi-directed graph.  A {\em minimally oriented graph\/} of
  $G$ is a graph $G^{\min}$ that satisfies the following three properties:
  \begin{enumerate}
  \item $G^{\min}$ is a maximal ancestral graph;
  \item $G$ and $G^{\min}$ are Markov equivalent;
  \item $G^{\min}$ has the minimum number of arrowheads of all
    maximal ancestral graphs that are Markov equivalent to $G$.  Here
    the number of arrowheads of an ancestral graph $G$ with $d$
    directed and $b$ bi-directed edges is defined as $\arr(G)=d+2b$. 
  \end{enumerate}
\end{definition}

By Lemma~\ref{lem:MAGs}, a minimally oriented graph $G^{\min}$ has the same
skeleton as the underlying bi-directed graph $G$.  According to
Theorem~\ref{thm:get-mconnect}, $G^{\min}$ has the boundary containment
property.  Examples of minimally oriented graphs are shown in Figure
\ref{fig:Gs}.  Given the small number of vertices of these graphs the claim
that these graphs are indeed minimally oriented graphs can be verified
directly.  The example of graph $G_1$ in Figure \ref{fig:Gs} also
illustrates that minimally oriented graphs are not unique.  By symmetry,
reversing the direction of the edge $v\to w$ in the depicted $G_1^{\min}$
yields a second minimally oriented graph for $G_1$.

We now turn to the problem of how to construct a minimally oriented graph.
Define a relation on the vertex set $V$ of the given bi-directed graph $G$
by letting $v \preccurlyeq_B w$ if $v=w$ or if $\Bd (v) \subsetneq \Bd (w)$
in $G$.  The relation $\preccurlyeq_B$ is a partial order and can thus be
extended to a total order $\le$ on $V$ such that the strict boundary
containment $\Bd(v) \subsetneq \Bd(w)$ implies that $v < w$.  In general,
the choice of such an extension to a total order is not unique.


\newtheorem{algorithm}[theorem]{Algorithm}
\begin{algorithm}
  \label{alg:construct-Gos}
  Let $G$ be a bi-directed graph, and $\le$ a total order on $V$ that
  extends the partial order $\preccurlyeq_B$ obtained from strict boundary
  containment.  Create a new graph $G^{\min}_<$ as follows:
  \begin{enumerate}
  \item[(a)] find the simplicial graph $G^s$ of $G$;
  \item[(b)] set $G^{\min}_<=G^s$;
  \item[(c)] replace every bi-directed edge $v\bi w\in G^{\min}_<$ with
    $\Bd(v)\subseteq\Bd(w)$ and $v<w$ by the directed edge $v\to
    w$.
  \end{enumerate}
\end{algorithm}

The notation $G^{\min}_<$ indicates the dependence of this graph on {\em
  both\/} the bi-directed graph $G$ and the total order $\leq$.  Clearly,
by Theorem~\ref{thm:get-mconnect}, in order for $G^{\min}_<$ to be a
minimally oriented graph it is necessary that it satisfies the boundary
containment property.  The next lemma shows that this is true.

\begin{lemma}  
  \label{lem:edgesinGos}
  Let $G$ be a bi-directed graph and $G^{\min}_<$ the graph constructed in
  Algorithm~\ref{alg:construct-Gos}.  It then holds that
  \begin{enumerate}
  \item if $v - w$ is an undirected edge in $G^{\min}_<$, then $\Bd(v)=\Bd(w)$;
  \item if $v\to w$ is a directed edge in $G^{\min}_<$, then
    $\Bd(v)\subseteq \Bd(w)$; 
  \item $v\bi w$ is a bi-directed edge in $G^{\min}_<$ iff
    $\Bd(v)\not\subseteq\Bd(w)\not\subseteq\Bd(v)$.  
  \end{enumerate} 
\end{lemma}
\begin{proof}
  (i) follows directly from Lemma~\ref{lem:orientedgesinGs}(i) because it
  follows from Algorithm \ref{alg:construct-Gos} that $G^{\min}_<$ and
  $G^s$ contain the same undirected edges.
  
  (ii) If the edge $v\to w$ is already present in $G^s$, then
  $\Bd(v)\subsetneq \Bd(w)$ according to
  Lemma~\ref{lem:orientedgesinGs}(ii).  If $v\to w$ is not already present
  in $G^s$, then $v<w$ and $\Bd(v)\subseteq \Bd(w)$.
  
  (iii) Suppose $v$ and $w$ are two adjacent vertices such that
  $\Bd(v)\not\subseteq\Bd(w)\not\subseteq\Bd(v)$.  Then $v\bi w$ in $G^{s}$
  and this edge cannot be replaced by a directed edge in step (c) of
  Algorithm~\ref{alg:construct-Gos}.  For the reversed claim, consider two
  adjacent vertices $v$ and $w$ such that $\Bd(v)\subseteq\Bd(w)$. (The
  other case is symmetric.)  If $v<w$, then according to the definition of
  the simplicial graph and step (c) of Algorithm~\ref{alg:construct-Gos}
  the edge between $v$ and $w$ in $G^{\min}_<$ cannot have an arrowhead at
  $v$ and thus cannot be bi-directed.  If $v>w$, then $\Bd(v)=\Bd(w)$
  because $\Bd(v)\subsetneq\Bd(w)$ would imply $v<w$.  It follows that the
  edge between $v$ and $w$ in $G^{\min}_<$ cannot be bi-directed as
  otherwise the arrowhead at $w$ would be removed in step (c).
\end{proof}

By Lemma~\ref{lem:edgesinGos}(iii), $v\bi w\in G^{\min}_<$ iff there exist
vertices $x\in\bd(v)\setminus\{w\}$ and $y\in\bd(w)\setminus\{v\}$ such
that the induced subgraph $G_{\{x,y,v,w\}}$ equals one of the two graphs
shown in Figure \ref{fig:nonorient}.  Graphs that do not contain the
four-cycle from Figure~\ref{fig:nonorient}(ii) as an induced subgraph are
known as chordal or decomposable and play an important role in graphical
modelling \citep{lau:bk}.  Graphs not containing the path from
Figure~\ref{fig:nonorient}(i) as an induced subgraph are called cographs
and have favorable computational properties \citep{brandstadt:1999}.  For
instance, cographs can be recognized in linear time \citep{corneil:1985}.

\begin{theorem}  
  \label{thm:equivGos}
  The graph $G^{\min}_<$ constructed in Algorithm~\ref{alg:construct-Gos} is a
  minimally oriented graph for the bi-directed graph $G$.  
\end{theorem}

\begin{proof}  We verify the conditions (i) and (iii) of
  Definition~\ref{def:Gos}.  This is sufficient because $G^{\min}_<$ has
  the boundary containment property (Lemma~\ref{lem:edgesinGos}) and thus
  condition (i) implies condition (ii) by Theorem~\ref{thm:get-mconnect}.
    \medskip

  \noindent {\em (i) $G^{\min}_<$ is a maximal ancestral graph:}\\
  By Lemma~\ref{lem:MAGs} it suffices to show that $G^{\min}_<$ is an
  ancestral graph.  Let $v$ and $w$ be adjacent vertices such that $v-w\in
  G^{\min}_<$.  This is equivalent to $v-w\in G^s$, and it follows that
  there does not exist an arrowhead at $v$ or $w$; compare the proof of
  Theorem \ref{thm:equivsimple}.  Furthermore, $G^{\min}_<$ does not
  contain any directed cycles because Algorithm~\ref{alg:construct-Gos}
  ensures that the presence of a directed edge $v\to w\in G^{\min}_<$
  implies $v<w$ in the total order.  Finally, assume that there exists
  $v\bi w\in G^{\min}_<$.  Then there cannot be a directed path from $v$ to
  $w$, since by Lemma \ref{lem:edgesinGos}(ii) this would imply
  $\Bd(v)\subseteq\Bd(w)$, contradicting Lemma \ref{lem:edgesinGos}(iii).
  \medskip

\begin{figure}[tbp]
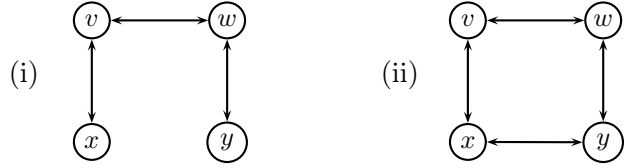

  \centering
  \setlength{\unitlength}{0.5pt}
  \begin{tabular}{c@{\hspace{1cm}}c}
    \parbox{4cm}{
      \vspace{2.75cm}   
      \small\psset{unit=1.2cm}
      \settowidth{\MyLength}{$v$}
      \newcommand{\myNode}[2]{\circlenode{#1}{\makebox[\MyLength]{#2}}}
      \rput(0.25,1.25){\normalsize (i)}
      \rput(1,1.85){\myNode{v}{$v$}} 
      \rput(1,0.5){\myNode{x}{$x$}} 
      \rput(2.5,1.85){\myNode{w}{$w$}}
      \rput(2.5,0.5){\myNode{y}{$y$}}
      \ncline{<->}{v}{w}
      \ncline{<->}{w}{y}
      \ncline{<->}{v}{x}
      \vspace{-0.1cm}
    }
    &
    \parbox{4cm}{
      \vspace{2.75cm}   
      \small\psset{unit=1.2cm}
      \settowidth{\MyLength}{$v$}
      \newcommand{\myNode}[2]{\circlenode{#1}{\makebox[\MyLength]{#2}}}
      \rput(0.25,1.25){\normalsize (ii)}
      \rput(1,1.85){\myNode{v}{$v$}} 
      \rput(1,0.5){\myNode{x}{$x$}} 
      \rput(2.5,1.85){\myNode{w}{$w$}}
      \rput(2.5,0.5){\myNode{y}{$y$}}
      \ncline{<->}{v}{w}
      \ncline{<->}{x}{y}
      \ncline{<->}{w}{y}
      \ncline{<->}{v}{x}
      \vspace{-0.1cm}
    }    
  \end{tabular}
  \caption{Induced subgraphs for which no
      arrowhead can be dropped from edge $v\bi w$.}
  \label{fig:nonorient}
\end{figure}

\noindent {\em (iii) $G^{\min}_<$ has the minimal number of arrowheads:}\\
Let $\bar G$ be a maximal ancestral graph that is Markov equivalent to the
(bi-directed) graph $G$, which requires that $\bar G$ and $G$, and thus
also $G^{\min}_<$ have the same skeleton.  Assume that $\arr(\bar G)<
\arr(G^{\min}_<)$.  Then either (a) there exists $v\to w\in G^{\min}_<$
such that $v-w\in \bar G$ or (b) there exists $v\bi w\in G^{\min}_<$ such
that $v\to w\in \bar G$ or $v-w\in \bar G$.
  
  Case (a): If $v\to w\in G^{\min}_<$, then $w$ cannot be simplicial.  Hence,
  there exist two vertices $x,y\in \bd(w)$ that are not adjacent in
  $G^{\min}_<$, and thus not adjacent in $G$; ($v=x$ is possible).  The global
  Markov property of $G$ states that $x\ind y$.  Since $\bar G$ is an
  ancestral graph and $v-w\in \bar G$, however, there may not be any
  arrowheads at $w$ on the edges between $x$ and $w$, and $y$ and $w$ in
  $\bar G$.  Therefore, $x$ and $y$ are $m$-connected given $\emptyset$ in
  $\bar G$, which yields that the global Markov property of $\bar G$ does
  not imply $x\ind y$; a contradiction.
  
  Case (b): Suppose $v\bi w\in G^{\min}_<$ but there is no arrowhead at $v$ on the
  edge between $v$ and $w$ in $\bar G$.  By Lemma \ref{lem:edgesinGos}(iii)
  there exists $x\in \bd(v)\setminus\Bd(w)$ such that $x$ and $w$ are not
  adjacent in $G^{\min}_<$.  Thus $x$ and $w$ are not adjacent in $G$ and
  $x\ind w$ is stated by the global Markov property for $G$.  In $\bar G$,
  however, $v$ is a non-collider on the path $(x,v,w)$ and thus this path
  $m$-connects $x$ and $w$ given $\emptyset$, which yields that the global
  Markov property of $\bar G$ does not imply $x\ind w$; a contradiction.
\end{proof}

The next result shows that our construction of minimally oriented graphs is
complete in the sense that every minimally oriented graph can be obtained
as the output of Algorithm~\ref{alg:construct-Gos} by appropriate choice of
a total order on the vertex set.

\begin{theorem}
  \label{thm:every-Gos}
  If $G^{\min}$ is a minimally oriented graph for a bi-directed graph $G$,
  then there exists a total order $\le$ on the vertex set such that
  $G^{\min}=G^{\min}_<$.
\end{theorem}
\begin{proof}
  The graph $G^{\min}$ is an ancestral graph and thus contains no directed
  cycles.  Hence, the directed edges in $G^{\min}$ yield a partial order
  $\preccurlyeq_D$ on the vertex set $V$ in which $v\preccurlyeq_D w$ if
  $v=w$ or if there is a directed path from $v$ to $w$.  Define the
  relation $\preccurlyeq_{\mathit{BD}}$ by letting $v\preccurlyeq_{\mathit{BD}} w$ if
  $v\preccurlyeq_{B} w$ or $v\preccurlyeq_{D} w$.  Clearly,
  $v\preccurlyeq_{\mathit{BD}} v$, i.e., the relation is reflexive.  We claim that
  the relation is in fact a partial order.
  
  By Theorem~\ref{thm:get-mconnect}, $G^{\min}$ has the boundary
  containment property such that $\Bd(v)\subseteq\Bd(w)$ if
  $v\preccurlyeq_D w$.  Consequently, if $v\not=w$ then $v\preccurlyeq_D w$
  implies $w\not\preccurlyeq_B v$ and $v\preccurlyeq_B w$ implies
  $w\not\preccurlyeq_D v$.  This implies that $\preccurlyeq_{\mathit{BD}}$ is
  anti-symmetric.  In order to verify transitivity, it suffices to consider
  three distinct vertices satisfying $v\preccurlyeq_{D} w\preccurlyeq_{B}
  u$ or $v\preccurlyeq_{B} w\preccurlyeq_{D} u$.  In the former case
  $\Bd(v)\subseteq\Bd(w)\subsetneq\Bd(u)$, and in the latter case
  $\Bd(v)\subsetneq\Bd(w)\subseteq\Bd(u)$.  In both cases
  $\Bd(v)\subsetneq\Bd(u)$ such that $v\preccurlyeq_B u$, which implies the
  required conclusion $v\preccurlyeq_{\mathit{BD}} u$.  
  
  We can now choose a total order $\le$ on $V$ that extends the partial
  order $\preccurlyeq_{\mathit{BD}}$ and thus extends both $\preccurlyeq_B$
  and $\preccurlyeq_D$.  Let $G^{\min}_<$ be the output of
  Algorithm~\ref{alg:construct-Gos} when the bi-directed graph $G$ and the
  chosen total order $\le$ are given as the input.  We claim that
  $G^{\min}=G^{\min}_<$.
  
  First note that if $v$ is a simplicial vertex of $G$, then there are no
  arrowheads at $v$ in $G^{\min}$.  Otherwise, we could drop all arrowheads
  at simplicial vertices in $G^{\min}$ to obtain an ancestral graph
  (Lemma~\ref{lem:drop-at-simplicial}) with fewer arrowheads.  The new
  graph would have the boundary containment property and thus be Markov
  equivalent to $G$ (by Theorem~\ref{thm:get-mconnect}).  This would
  contradict the assumed minimality of $G^{\min}$.
  
  The observation about simplicial vertices implies that an undirected edge
  in the simplicial graph $G^s$ is also an undirected edge in $G^{\min}$.
  Conversely, if $v-w\in G^{\min}$ then there may not be an arrowhead at
  $v$ on any other edge, and likewise for $w$, because $G^{\min}$ is
  ancestral.  Since $G^{\min}$ has the boundary containment property, it
  follows from Proposition~\ref{prop:completebd} that both $v$ and $w$ are
  simplicial vertices.  This implies that $v-w\in G^s$ and we conclude that
  $G^{\min}$ and $G^s$ have the same undirected edges.  By construction,
  the same holds for $G^{\min}_<$ and $G^s$.  Hence, $G^{\min}$ and
  $G^{\min}_<$ have the same undirected edges.
  
  Suppose $v\to w\in G^{\min}$.  Then $\Bd(v)\subseteq\Bd(w)$ because
  $G^{\min}$ has the boundary containment property.  Moreover, $v<w$
  because the total order $\le$ extends $\preccurlyeq_D$.  It follows that
  $v\to w\in G^{\min}_<$.  In other words, every directed edge in
  $G^{\min}$ is also in $G^{\min}_<$.  This together with the fact that
  $G^{\min}$ and $G^{\min}_<$ have the same skeleton and the same number of
  arrowheads, $\arr(G^{\min})=\arr(G^{\min}_<)$, implies that
  $G^{\min}=G^{\min}_<$.
\end{proof}

\subsection{Markov equivalence results}
\label{subsec:mini-orient-Mequiv}

The following corollary is an immediate consequence of Proposition
\ref{prop:equivsimple} because a minimally oriented graph $G^{\min}$ is an
undirected graph iff $G^s$ is an undirected graph.  

\begin{corollary}
  \label{cor:MarkovGos}
  Let $G^{\min}$ be a minimally oriented graph for a bi-directed graph
  $G$.  If $G$ is Markov equivalent to an undirected graph $U$, then
  $G^{\min}=U$ is the unique minimally oriented graph of $G$. 
\end{corollary}

A minimally oriented graph also reveals whether the original bi-directed
graph is Markov equivalent to a DAG.

\begin{theorem}
  \label{thm:MarkovGos}
  Let $G^{\min}$ be a minimally oriented graph for a bi-directed graph
  $G$.  Then $G$ is Markov equivalent to a DAG iff $G^{\min}$ contains no
  bi-directed edges.
\end{theorem}

\begin{proof}
  Let $G$ be a bi-directed graph such that $G^{\min}$ contains no
  bi-directed edges.  If $A\subseteq V$ is a simplicial set, then the
  induced subgraph $(G^{\min})_A$ is undirected and complete (this follows
  directly from Theorem \ref{thm:every-Gos} and Algorithm
  \ref{alg:construct-Gos}).  Let $A_1,\ldots, A_q$ be the inclusion-maximal
  simplicial sets of $G$.  Let $D$ be a directed graph obtained by
  replacing each induced subgraph $(G^{\min})_{A_i}$, $i=1,\ldots,q$, by a
  complete DAG.  Then $D$ itself has to be acyclic, which can be seen as
  follows: First, since $G^{\min}$ does not contain any directed cycles, a
  directed cycle $\pi$ in $D$ must involve a vertex $v\in \cup_{i=1}^q
  A_i$.  Let $v\in A_j$.  Since the induced subgraphs $D_{A_i}$,
  $i=1,\ldots,q$, are all acyclic, $\pi$ must also involve a vertex not in
  $A_j$. Therefore, there exists an edge $x\to w$ on $\pi$ such that $w\in
  A_j$ and $x\not\in A_j$.  Since the sets $A_i$ are inclusion-maximal
  simplicial sets, no vertex in $A_i$, $i\not= j$, is adjacent to any
  vertex in $A_j$.  Hence, $x\not\in\cup_{i=1}^q A_i$, which implies that
  the edge $x\to w$ is also present in $G^{\min}$.  This is a contradiction
  to $w$ being a simplicial vertex.
  
  Two vertices are adjacent in $G^{\min}$ iff they are adjacent in $D$.
  Moreover, $D$ has the boundary containment property because $G^{\min}$
  has this property, and if $u\rightarrow \bar{u}$ in $D$ then either $u
  \rightarrow \bar{u}$ in $G^{\min}$ or $u-\bar{u}$ in $G^{\min}$.  It
  thus follows from Theorem~\ref{thm:get-mconnect} that $D$ is Markov
  equivalent to $G^{\min}$ and $G$.
  
  Conversely, suppose that $v\bi w\in G^{\min}$ and for a contradiction,
  that $G$ is Markov equivalent to a DAG $D$.  Note that $D$ must have the same
  skeleton as $G$ (and $G^{\min}$). By Lemma \ref{lem:edgesinGos}(iii),
  there exist two different vertices $x\in\bd(v)\setminus\{w\}$ and
  $y\in\bd(w)\setminus\{v\}$ such that, by the Markov property of $G$,
  $x\ind w$ and $v\ind y$.  Hence, $v$ and $w$ must be colliders on the
  paths $(x,v,w)$ and $(v,w,y)$ in $D$, respectively.  This is impossible
  in the DAG $D$.
\end{proof}

Theorem \ref{thm:MarkovGos} can be shown to be equivalent to a Markov
equivalence result stated without proof in Theorem~1 in \cite{pearl:1994}.
This latter theorem requires `no chordless four-chain', which must be read
as excluding graphs with induced subgraphs that are either of the graphs in
Figure~\ref{fig:nonorient}.  Under this condition, \cite{pearl:1994} also
state that a Markov equivalent DAG can be constructed from the (undirected)
skeleton of $G$ by introducing directed and bi-directed edges in an
operation they term `sink orientation', and turning remaining undirected
edges into directed ones.  The sink orientation of the graph $G_1$ in
Figure~\ref{fig:Gs} has the directed edges of $G_1^s$ but an undirected
edge $v-w$. Thus sink orientation need not yield an ancestral graph.
The bi-directed graphical models considered in Theorem \ref{thm:MarkovGos}
also appear in the construction of generalized Wishart distributions
\citep[][Thm.~2.2]{letac:massam:wishart:2005}.  In that context the models
are called {\it homogeneous} and characterized in terms of Hasse diagrams.

As the next result reveals, bi-directed graphs that are Markov equivalent
to DAGs exhibit a structure that corresponds to a multivariate regression
model.  The graphs can also be termed chordal cographs; compare the
paragraph before Theorem~\ref{thm:equivGos}.

\begin{proposition}
  Let $G^{\min}$ be a minimally oriented graph for a connected bi-directed
  graph $G$.  If $G^{\min}$ contains no bi-directed edges, then the set $A$
  of all simplicial vertices is non-empty, the induced subgraph
  $(G^{\min})_A$ is a disjoint union of complete undirected graphs, the
  induced subgraph $(G^{\min})_{V\setminus A}$ is a complete DAG, and an
  edge $v\to w$ joins any two vertices $v\in A$ and $w\not\in A$ in
  $G^{\min}$.
\end{proposition}
\begin{proof}
  For two adjacent vertices $v$ and $w$ in $G^{\min}$, Lemma
  \ref{lem:edgesinGos}(i)-(ii) implies that $\Bd(v)\subseteq\Bd(w)$ or
  $\Bd(w)\subseteq\Bd(v)$.  Hence, we can list the vertex set as
  $V=\{v_1,\dots,v_p\}$ such that $\Bd(v_i)\subseteq\Bd(v_j)$ if $v_i$ and
  $v_j$ are adjacent and $i\le j$.  It follows that $v_1\in A$ and thus
  $A\not=\emptyset$.  Let $A_1,\ldots, A_q$ be the inclusion-maximal
  simplicial sets of $G$.  Then $(G^{\min})_A$ equals the union of the
  disjoint complete undirected graphs
  $(G^{\min})_{A_1},\dots,(G^{\min})_{A_q}$.  Since $G^{\min}$ is an ancestral
  graph, $(G^{\min})_{V\setminus A}$ is a DAG.
  
  We prove the remaining claims by induction on $|V\setminus A|$.  If
  $|V\setminus A|=0$, then the connected graph $G^{\min}$ is a complete
  undirected graph and there is nothing to show.  Let $|V\setminus A|\ge
  1$.  It follows that $v_p \in V\setminus A$.  If the shortest path
  between some vertex $v_{i_1}$ and $v_p$ in $G$ is of the form $v_{i_1}\bi
  \dots\bi v_{i_k}\bi v_p$, then $i_1<\dots<i_k<p$ and
  $\Bd(v_{i_1})\subseteq\dots \subseteq\Bd(v_{i_k})\subseteq\Bd(v_p)$,
  which is easily shown by induction on $k$.  However, since
  $v_{i_1}\in\Bd(v_{i_1})$ it must in fact hold that $v_{i_1}$ and $v_p$
  are adjacent.  Hence, there is an edge between every vertex $v\in
  V\setminus \{v_p\}$ and $v_p$, which for $v\in A$ is of the form $v\to
  v_p$ because clearly $v_p\not\in A$. The proof is finished by combining
  what we learned about $v_p$ with the induction assumption applied to the
  induced subgraph $G_W$ with $W={\{v_1,\dots,v_{p-1}\}}$.  Note that for
  $v,w\in W$, the inclusion $\Bd_G(v)\subseteq\Bd_G(w)$ implies that
  $\Bd_{G_W}(v)\subseteq\Bd_{G_W}(w)$. Thus by Lemma \ref{lem:edgesinGos}
  and Theorem \ref{thm:equivGos}, $(G_W)^{\min}$ does not contain any
  bi-directed edges.
\end{proof}

\section{Maximum likelihood estimation in Gaussian models} 
\label{sec:impl-max-lik}

In this section we consider the Gaussian covariance models associated with
bi-directed graphs and demonstrate that the graphical constructions from
Sections~\ref{sec:simpl-sets} and \ref{sec:orient-edges} can be employed
for more efficient computation of maximum likelihood estimates.

\subsection{Covariance graphs and Gaussian ancestral graph models}
\label{subsec:covgraphmodels}

Let $G$ be a bi-directed graph, and
\begin{equation}
  \label{eq:PG}
  \mathbf{P}(G) = \bigl\{ \Sigma\in \RRR^{V\times V} \mid
  \Sigma=(\sigma_{vw}) \mbox{ sym.\ pos.\  def.},\; \sigma_{vw}=0 \;\forall
  (v,w): v\bi w\notin G\bigr\} 
\end{equation}
be the cone of symmetric positive definite matrices with zero pattern
induced by $G$.  The {\em covariance graph model\/} associated with $G$ is
the family of multivariate normal distributions $\mathbf{N}(G) = \bigl(
\ND_V(0,\Sigma) \mid \Sigma\in \mathbf{P}(G)\bigr)$.  It can be shown that
every distribution in $\mathbf{N}(G)$  satisfies all conditional
independences stated by the global Markov property for the bi-directed
graph $G$ \citep[Prop.\ 2.2]{kauermann:dual}.  Conversely, if a
distribution $\ND_V(0,\Sigma)$ satisfies the global Markov property for
$G$, then $\Sigma\in\mathbf{P}(G)$.

Let $S\in\RRR^{V\times V}$ be the empirical covariance matrix computed from
an i.i.d.\ sample drawn from some unknown distribution
$\ND_V(0,\Sigma)\in\mathbf{N}(G)$, i.e., the $(v,w)$-th entry in $S$ is the
dot-product of the vectors of observations for the $v$-th and $w$-th
variables divided by the sample size $n$.  The log-likelihood function
$\ell_{S,n} : \mathbf{P}(G) \to \RRR$ of $\mathbf{N}(G)$ can be written as
\begin{equation}
  \label{eq:loglik}
      \ell_{S,n}(\Sigma) = 
      -\frac{n|V|}{2}\log(2\pi) -\frac{n}{2}\log |\Sigma| -
      \frac{n}{2}\tr\left(\Sigma^{-1}S\right).
\end{equation}
If $S$ is positive definite then the global maximum of $\ell_{S,n}$ over
$\mathbf{P}(G)$ exists. 
The likelihood equations obtained by setting to zero the partial
derivatives of $\ell_{S,n}$ with respect to the non-restricted entries in
$\Sigma$ take on the form
\begin{equation}\label{eq:likeqn}
\begin{split}
(\Sigma^{-1})_{vw}&=(\Sigma^{-1}S\Sigma^{-1})_{vw} \quad \forall
v,w\in V: \; v=w\;\; \mathrm{ or }\;\; v\bi w\in G;
\end{split}
\end{equation}
compare \citet[\S 2.1.1]{anderson:1985}.  A matrix
$\hat\Sigma(S)\in\mathbf{P}(G)$ that solves (\ref{eq:likeqn}) is {\em a
  solution to the likelihood equations\/} of $\mathbf{N}(G)$.  Since
subsequent theorems on the structure of the likelihood equations are
obtained via Gaussian ancestral graph models, we briefly review the
parametrization of these models.


Let $G$ be an ancestral graph and $\un_G\subseteq V$ the set of vertices
$v$ that are such that any edge with endpoint $v$ has a tail at $v$.  By
Definition \ref{def:ancestralgraph}(i), $v-w\in G$ implies $v,w\in \un_G$,
and $v\bi w\in G$ implies that $v,w\notin \un_G$.  Let $\Lambda$ be a
symmetric positive definite $\un_G\times \un_G$ matrix such that
$\Lambda_{vw}\not=0$ only if $v= w$ or $v-w\in G$.  Let $\Omega$ be a
symmetric positive definite $(V\setminus\un_G)\times (V\setminus \un_G)$
matrix such that $\Omega_{vw}\not=0$ only if $v= w$ or $v\bi w\in G$.
Finally, let $B$ be a $V\times V$ matrix such that $B_{vw}\not= 0$ only if
$w\to v\in G$.  Define the symmetric positive definite matrix
\begin{equation}
  \label{eq:sigma}
  \Sigma(\Lambda,B,\Omega) = (I-B)^{-1}\begin{pmatrix} \Lambda^{-1}
  & 0 \\ 0 & 
  \Omega\end{pmatrix} (I-B)^{-{\rm T}},
\end{equation}
where $I$ is the identity matrix.

Let $\mathbf{N}(G)$ be the Gaussian ancestral graph model associated with
$G$, i.e., the family of all centered normal distributions that are
globally Markov with respect to $G$.  As shown in
\citet[\S8]{richardson:2002}, the normal distribution $\ND_V(0,\Sigma)$
with $\Sigma=\Sigma(\Lambda,B,\Omega)$ defined in (\ref{eq:sigma}) is in
$\mathbf{N}(G)$.  Conversely, if $G$ is {\em maximal\/},
then for any $\ND_V(0,\Sigma)\in\mathbf{N}(G)$ there exist unique $\Lambda,
\Omega, B$ of the above type such that $\Sigma=\Sigma(\Lambda,B,\Omega)$.
(Note that \cite{richardson:2002} use $B$ for what is here denoted
by $I-B$.)

Since a bi-directed graph $G$ and a minimally oriented graph $G^{\min}$ are
Markov equivalent, the parametrization map for $G^{\min}$,
$(\Lambda,B,\Omega) \mapsto\Sigma(\Lambda,B,\Omega)$, has image equal to
$\mathbf{P}(G)$.  By \citet[Thm.\ 8.14, Lemma 8.22]{richardson:2002}, we
obtain the following Lemma.
\begin{lemma}
  \label{lem:conserve-stationary}
  Let $G$ be a bi-directed graph.  The covariance matrix
  $\Sigma(\Lambda,B,\Omega)$ solves the likelihood equations of
  $\mathbf{N}(G)$ iff $(\Lambda,B,\Omega)$ solves the likelihood equations
  of $\mathbf{N}(G^{\min})$.
\end{lemma}

\subsection{Empirical maximum likelihood estimates}

Using the graphical results established earlier, we can show that over
simplicial sets a solution to the likelihood equations (\ref{eq:likeqn})
agrees with its empirical counterpart in $S$.

\begin{theorem}  
  \label{thm:mle}
  Let $G$ be a bi-directed graph with associated covariance graph
  model $\mathbf{N}(G)$. If $A\subseteq V$ is simplicial, $S$
  is a symmetric positive definite matrix, and 
  $\hat\Sigma(S)\in\mathbf{P}(G)$ is a solution to the likelihood
  equations (\ref{eq:likeqn}),  then
$\hat\Sigma(S)_{A\times A} = S_{A\times A}$.
\end{theorem}
\begin{proof}
  By Theorem \ref{thm:equivsimple}, the covariance graph model
  $\mathbf{N}(G)$ and the Gaussian ancestral graph model
  $\mathbf{N}(G^{s})$ based on the simplicial graph $G^{s}$ are equal.  Let
  $\mathbf{N}(G^{s})$ be parametrized by the precision matrix $\Lambda$,
  the matrix of regression coefficients $B$ and the covariance matrix
  $\Omega$ as described in \S \ref{subsec:covgraphmodels}.  In particular,
  it follows from \citet[Lemma 8.4]{richardson:2002} that if
  $\Sigma=\Sigma(\Lambda,B,\Omega)$, then $(\Lambda^{-1})_{A\times
    A}=\Sigma_{A\times A}$.
  
   The inclusion-maximal simplicial sets $A_1,\ldots,A_q$ of $G$ form a
   partition of $\un_{G^s}$.  The induced subgraphs $G^{s}_{A_i}$,
   $i=1,\ldots,q$, are complete undirected graphs.  It follows that
   $\Lambda$ is a block-diagonal matrix such that $\Lambda_{vw}=0$ if there
   does not exist an inclusion-maximal simplicial set $A_i$ such that
   $v,w\in A_i$.  Now the discussion in \citet[\S8.5]{richardson:2002} and
   Lemma \ref{lem:conserve-stationary} imply that every solution to the
   likelihood equations for $\Lambda$, $B$, $\Omega$ in the Gaussian
   ancestral graph model $\mathbf{N}(G^{s})$ satisfies that
   $(\hat{\Lambda}^{-1})_{A_i\times A_i} = S_{A_i\times A_i}$ for all
   $i=1,\ldots,q$.  Since $A\subseteq A_j$ for some $j$, it holds that
   $\hat\Sigma_{A\times A}=(\hat{\Lambda}^{-1})_{A\times A} = S_{A\times
     A}$.
\end{proof}


Our graphical constructions also provide information on when maximum
likelihood estimates of conditional parameters are equal to their empirical
counterparts. The conditional parameters we consider are the regression
coefficients and conditional variance for the conditional distribution of
variable $v$ given its {\em parents\/} $ \pa(v)=\{w\in V\mid w\to v\in
G^{\min}\}$ in a minimally oriented graph $G^{\min}$.  If $\pa(v)=\emptyset$,
then conditioning variable $v$ on $\pa(v)$ is understood to yield the
marginal distribution of $v$.
\begin{theorem}  
  \label{thm:mle3}
  Let $G^{\min}$ be a minimally oriented graph for a bi-directed graph
  $G$, $S$ a symmetric positive definite matrix, and
  $\hat\Sigma(S)\in\mathbf{P}(G)$ a solution to the likelihood equations
  (\ref{eq:likeqn}).  If $v$ is a vertex such that there is no vertex $w$
  with $v\bi w\in G^{\min}$, then the regression coefficients for $v$ given
  $\pa(v)$ are
  \begin{equation}\label{eq:mle3a}
  \hat\Sigma(S)_{v\times \pa(v)} 
  \bigl[\hat\Sigma(S)_{\pa(v)\times \pa(v)}\bigr]^{-1} 
  = S_{v\times \pa(v)}\bigl(S_{\pa(v)\times \pa(v)}\bigr)^{-1},
  \end{equation}
  and that the conditional variance for $v$ given
  $\pa(v)$ is
  \begin{multline}\label{eq:mle3b}
      \hat\Sigma(S)_{vv}-\hat\Sigma(S)_{v\times \pa(v)} 
      \bigl[\hat\Sigma(S)_{\pa(v)\times \pa(v)}\bigr]^{-1}
      \hat\Sigma(S)_{\pa(v)\times v}
      = \\
      S_{vv}-S_{v\times \pa(v)}\bigl(S_{\pa(v)\times \pa(v)}\bigr)^{-1}
      S_{\pa(v)\times v}.
  \end{multline}
\end{theorem}
\begin{proof}
  If $\pa(v)=\emptyset$, then $v$ is a simplicial vertex, and the claim
  reduces to $\hat\Sigma(S)_{vv}=S_{vv}$, which follows from Theorem
  \ref{thm:mle}.  Otherwise, using the parametrization of
  $\mathbf{N}(G^{\min})$, it follows from \citet[Thm.\ 8.7]{richardson:2002}
  that if $\Sigma=\Sigma(\Lambda,B,\Omega)$, then
  \[
  \Sigma_{v\times \pa(v)} \bigl[\Sigma_{\pa(v)\times
    \pa(v)}\bigr]^{-1}= B_{v\times \pa(v)}
  \]
  and
  \[
  \Sigma_{vv}-\Sigma_{v\times \pa(v)} 
  \bigl[\Sigma_{\pa(v)\times \pa(v)}\bigr]^{-1}
  \Sigma_{\pa(v)\times v} 
  = \Omega_{vv}.
  \]
  If $\hat{\Lambda}$, $\hat{B}$, $\hat{\Omega}$ solve the likelihood
  equations for $\mathbf{N}(G^{\min})$, then $\hat{B}_{v\times \pa(v)}$ and
  $\hat{\Omega}_{vv}$ solve the likelihood equations of the model in which
  all parameters in $\Lambda$, $B$, $\Omega$ except for $B_{v\times
    \pa(v)}$ and $\Omega_{vv}$ are held fixed.  It follows from
  \citet[\S\S5.1-2]{drton:2004c} that $\hat B_{v\times \pa(v)}$ and
  $\hat\Omega_{vv}$ are equal to the empirical expressions on the right
  hand side of (\ref{eq:mle3a}) and (\ref{eq:mle3b}), respectively.
  Applying Lemma \ref{lem:conserve-stationary} yields the claim.
\end{proof}

\begin{remark} 
  \label{rem:efficient}
  \rm {\em Iterative Conditional Fitting\/} is a special purpose algorithm
  for maximum likelihood estimation in covariance graph models
  \citep{drton:2003b,chaudhuri:2007}. However, it does not exploit the
  results of Theorems \ref{thm:mle} and \ref{thm:mle3}.  On the other hand,
  if one runs the ancestral graph extension of iterative conditional
  fitting described in \cite{drton:2004c} on a minimally oriented graph,
  then unnecessary computations are avoided by implicitly exploiting
  Theorems \ref{thm:mle} and \ref{thm:mle3}.  This is illustrated in the
  example in Section~\ref{sec:example}.
\end{remark}

If a bi-directed graph $G$ has a minimally oriented graph $G^{\min}$
without bi-directed edges then $G$ is Markov equivalent to a DAG
(Theorem~\ref{thm:MarkovGos}) and the likelihood equations have a unique
solution that is a rational function of the empirical covariance matrix
$S$.  However, this is no longer true if there is a bi-directed edge in
$G^{\min}$.  In this case, $G$ contains one of the two graphs in
Figure~\ref{fig:nonorient} as a subgraph; compare
Lemma~\ref{lem:edgesinGos}(iii).  Solving the likelihood equations for the
bi-directed four-chain in Figure~\ref{fig:nonorient}(i) is equivalent to
computing the roots of a quintic polynomial.  There exist data for which
this quintic has exactly three real roots \citep{drton:2004}.  Galois
theory \citep[Lemma 14.7]{stewart:1989} implies that for these data the
quintic is unsolvable by radicals, i.e., the roots of the quintic and thus
the solutions to the likelihood equations cannot be computed from the data
in finitely many steps involving addition, subtraction, multiplication,
division, or taking $r$-th roots.  (\cite{geiger:2006} obtain similar
results in the context of undirected graphs.) Similarly, solving the
likelihood equations of the bi-directed four-cycle in
Figure~\ref{fig:nonorient}(ii) corresponds to solving a polynomial equation
system of degree 17.  This can be verified in computer algebra systems such
as {\tt Singular} \citep{singular}; see also
\citet[\S5]{drtonsullivant:2007}.  It is natural to conjecture that there
exist data for which this system is also unsolvable by radicals.

\begin{table}[t]
  \centering
  
  \begin{tabular}{lrrrrrrrr}
    &  \mbox{\footnotesize GAL7} & \mbox{\footnotesize GAL10}  &
    \mbox{\footnotesize GAL1}  & \mbox{\footnotesize GAL3}  &
    \mbox{\footnotesize GAL2} & \mbox{\footnotesize GAL80} &
    \mbox{\footnotesize GAL11}  & \mbox{\footnotesize GAL4}\\
\hline\hline
\mbox{\footnotesize GAL7}  & \em 1.000 &0.91 &0.88 &0.50 &0.81 &0.21 &-0.07 &-0.08\\
\mbox{\footnotesize GAL10} &0.910 &\em 1.000 &0.92 &0.46 &0.87 &0.26 &-0.08 &-0.07\\
\mbox{\footnotesize GAL1}  &0.880 &0.920 &\em 1.000 &0.39 &0.87 &0.28 &-0.10 &-0.10\\
\mbox{\footnotesize GAL3}  &0.489 &0.447 &0.374 &\em 0.998 &0.44 &0.20 &-0.18 &0.12\\
\mbox{\footnotesize GAL2}  &0.807 &0.865 &0.865 &0.422 &\em 0.991 &0.26 &-0.18
    &-0.03\\
\mbox{\footnotesize GAL80} &0.224 &0.271 &0.297 &0.191 &0.280 &\em 1.001 &0.08 &0.23\\
\mbox{\footnotesize GAL11} &0 &0 &0 &-0.208 &-0.103 &0 &\em 1.022 &0.24\\
\mbox{\footnotesize GAL4}  &0 &0 &0 &0 &0.038 &0.209 &0.255 &\em 0.987\\
  \end{tabular}
  \caption{\label{tab:gene}
    Gene expression data. Empirical correlation matrix (above diagonal) and maximum
    likelihood estimate (below diagonal).  The italicized diagonal entries are ratios
    between  are maximum likelihood and empirical variance estimates. }
\end{table}

\subsection{Example: Gene expression measurements}
\label{sec:example}

The application of covariance graph models to gene expression data has been
promoted in \cite{butte:2000}.  For illustration, we
select data from microarray experiments with yeast strands
\citep{gasch:00}.  We focus on eight genes involved in galactose
utilization.  Expression measurements for all eight genes are available in
$n=134$ experiments, for which the empirical correlation matrix is shown in
the upper-diagonal part of Table \ref{tab:gene}.

\begin{figure}[t]
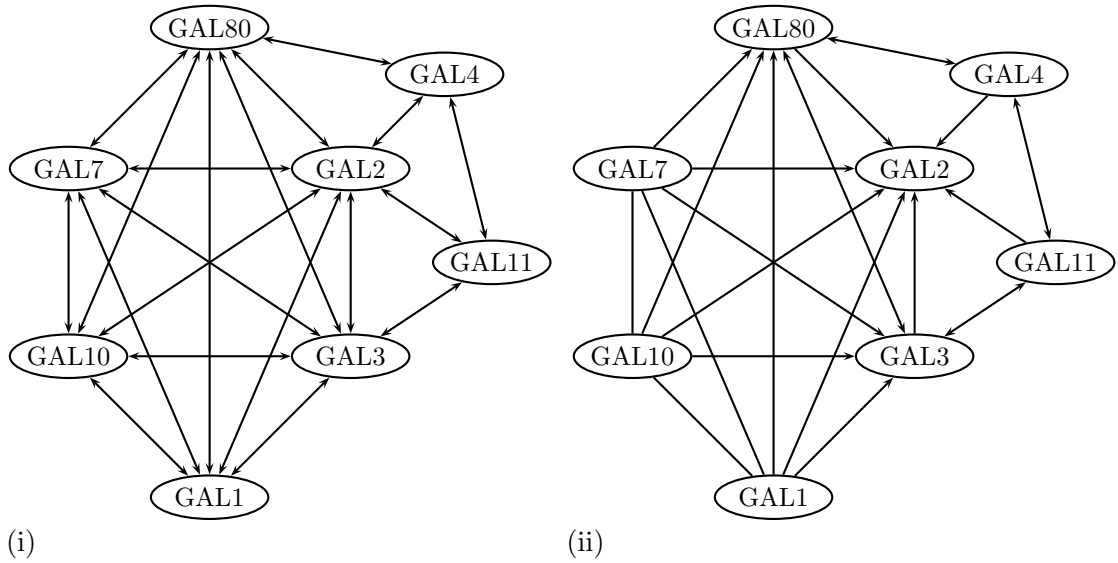

  \centering
  \begin{tabular}{l@{\hspace{1.5cm}}l}
    \hspace{-3.5cm}
      \parbox{6cm}{
      \vspace{7.5cm}   
      \small\psset{unit=1.25cm}
      \settowidth{\MyLength}{GAL1}
      \newcommand{\myNode}[2]{\ovalnode{#1}{\makebox[\MyLength]{#2}}}
      \rput(1,0){\normalsize (i)}
      \rput(1.5,4){\myNode{1}{GAL7}} 
      \rput(1.5,2){\myNode{2}{GAL10}} 
      \rput(3,0.5){\myNode{3}{GAL1}} 
      \rput(4.5,2){\myNode{4}{GAL3}} 
      \rput(4.5,4){\myNode{5}{GAL2}} 
      \rput(3,5.5){\myNode{6}{GAL80}} 
      \rput(6,3){\myNode{7}{GAL11}} 
      \rput(5.5,5){\myNode{8}{GAL4}} 
      \ncline{<->}{1}{2}
      \ncline{<->}{1}{3}
      \ncline{<->}{1}{4}
      \ncline{<->}{1}{5}
      \ncline{<->}{1}{6}
      \ncline{<->}{2}{3}
      \ncline{<->}{2}{4}
      \ncline{<->}{2}{5}
      \ncline{<->}{2}{6}
      \ncline{<->}{3}{4}
      \ncline{<->}{3}{5}
      \ncline{<->}{3}{6}
      \ncline{<->}{4}{5}
      \ncline{<->}{4}{6}
      \ncline{<->}{4}{7}
      \ncline{<->}{5}{6}
      \ncline{<->}{5}{7}
      \ncline{<->}{5}{8}
      \ncline{<->}{6}{8}
      \ncline{<->}{7}{8}
      \vspace{0.25cm}
    }
    &
      \parbox{6cm}{
      \vspace{7.5cm}   
      \small\psset{unit=1.25cm}
      \settowidth{\MyLength}{GAL1}
      \newcommand{\myNode}[2]{\ovalnode{#1}{\makebox[\MyLength]{#2}}}
      \rput(1,0){\normalsize (ii)}
      \rput(1.5,4){\myNode{1}{GAL7}} 
      \rput(1.5,2){\myNode{2}{GAL10}} 
      \rput(3,0.5){\myNode{3}{GAL1}} 
      \rput(4.5,2){\myNode{4}{GAL3}} 
      \rput(4.5,4){\myNode{5}{GAL2}} 
      \rput(3,5.5){\myNode{6}{GAL80}} 
      \rput(6,3){\myNode{7}{GAL11}} 
      \rput(5.5,5){\myNode{8}{GAL4}} 
      \ncline{-}{1}{2}
      \ncline{-}{1}{3}
      \ncline{->}{1}{4}
      \ncline{->}{1}{5}
      \ncline{->}{1}{6}
      \ncline{-}{2}{3}
      \ncline{->}{2}{4}
      \ncline{->}{2}{5}
      \ncline{->}{2}{6}
      \ncline{->}{3}{4}
      \ncline{->}{3}{5}
      \ncline{->}{3}{6}
      \ncline{->}{4}{5}
      \ncline{<->}{4}{6}
      \ncline{<->}{4}{7}
      \ncline{<-}{5}{6}
      \ncline{<-}{5}{7}
      \ncline{<-}{5}{8}
      \ncline{<->}{6}{8}
      \ncline{<->}{7}{8}
      \vspace{0.25cm}
    }
    \end{tabular}
  \caption{\label{fig:gene} (i) Bi-directed graph $G$ for gene expression
    measurements, (ii) the unique minimally oriented graph $G^{\min}$.}
\end{figure}

For these data, the covariance graph model induced by the graph $G$ in
Figure \ref{fig:gene}(i) has a deviance of 8.87 over 8 degrees of freedom,
which indicates a good model fit; the p-value computed using a chi-square
distribution is 0.35.  Figure \ref{fig:gene}(ii) shows the unique minimally
oriented graph $G^{\min}$.  The maximum likelihood estimate obtained by
fitting the model to the correlation matrix is shown in the lower-diagonal
part of Table \ref{tab:gene}; note that this estimate is not a correlation
matrix (not all the italicized diagonal entries are equal to one).  As
predicted by Theorem \ref{thm:mle}, the submatrix over GAL1, GAL7, and
GAL10 equals the respective submatrix in the empirical correlation matrix.
The regression coefficients for the regression of GAL2 on all remaining
variables are identical when computed from the maximum likelihood versus
the empirical estimate (Theorem \ref{thm:mle3}).

The use of a minimally oriented graph $G^{\min}$ leads to a considerable
gain in computational efficiency in the iterative calculation of the
maximum likelihood estimate $\hat\Sigma$.
With the identity matrix as starting value, iterative conditional fitting
(Remark~\ref{rem:efficient}) on the original bi-directed graph $G$ performs
eight multiple regressions per iteration and converges after 103
iterations.  Using the same starting value and termination criterion,
iterative conditional fitting on $G^{\min}$ converges after only 5
iterations and requires only five multiple regressions per iteration (for
the genes GAL2, GAL3, GAL4, GAL11, and {GAL80}), of which the one for GAL2
has to be executed only in the first iteration.

As in any application of covariance graph models, one might question the
assumption of Gaussianity.  Indeed there are $10$ experiments in which the
measurements for the genes GAL1, GAL7, GAL10 and GAL80 come out to be large
negative values, and one in which GAL7 alone takes such a value. These
appear to be outliers (standardized values between -3 and -5), possibly
produced by thresholding, as some values are identical. However, the
measurements for the other genes are well within the range of the
observations for the remaining 123 experiments. Thus it is unclear whether
removing these 11 experiments from consideration is appropriate.  If the 11
experiments are removed, then the correlations among GAL1, GAL7 and GAL10
decrease to values between $0.38$ and $0.60$, the latter value is the
maximum of all correlations.  Nevertheless the deviance for $G$ increases
only slightly to $10.09$ (p-value 0.26). The iterative conditional fitting
algorithm based on $G$ now converges after only $20$ iterations rather than
$103$. However, this is still four times as many iterations as required in
iterative conditional fitting based on the minimally oriented graph
$G^{\min}$; recall that in addition each iteration is also simpler.

The original correlation matrix in Table \ref{tab:gene} exhibits an apparent
similarity of the rows for GAL1, GAL7 and GAL10; this is also reflected in
the graph $G$ in which these variables form a complete set and have the
same spouses. Such symmetry could be investigated further via a group
symmetry model \citep{andersson:madsen:symmetry:1998}.

 \section{Conclusion}
\label{sec:discus}

We showed how to remove a maximal number of arrowheads from the edges of a
bi-directed graph $G$ such that one obtains a maximal ancestral graph
$G^{\min}$ that is Markov equivalent to $G$.  The graph $G^{\min}$, called
a minimally oriented graph, reveals whether $G$ is Markov equivalent to an
undirected graph, and also whether $G$ is Markov equivalent to a DAG.

For the (Gaussian) covariance graph model associated with $G$, a minimally
oriented graph $G^{\min}$ yields an alternative parametrization that
provides insight into likelihood inference.  The structure of the
arrowheads in $G^{\min}$ allowed us to identify parts of the covariance
matrix for which the maximum likelihood estimates are equal to their
empirical counterparts (this applies to all solutions to the likelihood
equations if, as occasionally happens, there is more than one solution).
This makes it possible to avoid or speed up iterative estimation of the
full covariance matrix.  We also saw that the maximum likelihood estimator
of the covariance matrix in a covariance graph model is a rational function
of empirical covariance matrix iff $G^{\min}$ contains no bi-directed edge.
This is similar to the results that identify decomposable models as the
sub-class of all log-linear and all covariance selection models
\citep{dempster:covsel} for which the maximum likelihood estimator is
available in closed form.

\citet{mdtsr:05} formulate binary models based on the Markov property of
bi-directed graphs.  For these models, the maximum likelihood estimator is
available in closed form if the model-inducing graph is Markov equivalent
to a DAG.  Moreover, we verified that in the example of the graph $G$ in
Figure \ref{fig:intro}, the maximum likelihood estimates of the marginal
distributions of $X_1$ and $X_4$ are equal to the corresponding empirical
proportions.  We thus believe that analogs to the Gaussian results
established here will hold in discrete models, but a general
parametrization of discrete ancestral graph models is required to fully
access the potential of the results obtained in this paper.

\begin{appendix}
\section{Connecting paths and boundary containment}
\label{sec:connect}

In this appendix we prove results about graphs that satisfy the {\em
  boundary containment property} from
Definition~\ref{def:dir-bound-contain}.  These results are used in the
proof of Theorem \ref{thm:get-mconnect}.

Let $v$ and $w$ be two fixed distinct vertices that are $m$-connected given
$C\subseteq V\setminus\{v,w\}$ in a simple mixed graph $G$.  Define
$\Pi_{G}(v,w|C)$ to be the set of paths that $m$-connect $v$ and $w$ given
$C$ in $G$, and let $\Pi_{G}^{\rm min}(v,w|C)$ be the set of paths that are
of minimal length among the paths in $\Pi_{G}(v,w|C)$.

\begin{lemma}\label{lem:connect1}
  If a simple mixed graph $G$ satisfies the boundary containment property,
  $v_{i-1}$, $v_i$ and $v_{i+1}$ are three consecutive vertices on a path
  $\pi$ in $G$, and $v_i$ is a non-collider on $\pi$, then $v_{i-1}$ and
  $v_{i+1}$ are adjacent.
\end{lemma}
\begin{proof}
  If $v_i$ is a non-collider, then the edge between $v_i$ and $v_{i-1}$ or
  the edge between $v_i$ and $v_{i+1}$ must have a tail at $v_i$.  Suppose,
  without loss of generality, that the latter is the case.  Then
  $\Bd(v_i)\subseteq\Bd(v_{i+1})$ and thus
  $v_{i-1}\in\Bd(v_{i+1})$, which is the claim.
\end{proof}

\begin{lemma}\label{lem:tsr1}
  Let $G$ be a simple mixed graph, and $\pi=(v,v_1,\ldots,v_k,w)\in
  \Pi_{G}(v,w|C)$.  Let $v_0=v$ and $v_{k+1}=w$.  If $v_i$ is a
  non-endpoint vertex on $\pi$ and there is an arrowhead at $v_i$ on the
  edge between $v_{i-1}$ and $v_i$, then either (i) $v_i \in \an(C)$ or
  (ii) the path $(v_i,v_{i+1},\ldots,v_k,w)$ is a directed path from $v_i$
  to $w$.
\end{lemma}
\begin{proof}
  Suppose the result is false. Let $v_j$ be the vertex closest to $w$
  satisfying the antecedent of the Lemma, but not the conclusion.  If $v_j$
  is a collider, then by definition of $m$-connection, $v_j \in \an (C)$,
  which is a contradiction.  If $v_j$ is a non-collider then $v_j \to
  v_{j+1}$ on $\pi$. If $v_{j+1}=w$, if $v_{j+1} \in \an (C)$, or if
  $(v_{j+1},\ldots,v_k,w)$ is a directed path from $v_{j+1}$ to $w$, then
  clearly $v_j$ satisfies the conclusion of the Lemma, which is a
  contradiction. But if $v_{j+1} \notin \an (C)\cup\{w\}$ and
  $(v_{j+1},\ldots,v_k,w)$ is not a directed path from $v_{j+1}$ to $w$
  then $v_{j+1}$ satisfies the conditions on $v_j$, but is closer to $w$,
  again a contradiction.
\end{proof}

\begin{lemma}\label{lem:tsr2}
  If $G$ is an ancestral graph that satisfies the boundary
  containment property and $\pi=(v,v_1,\ldots,v_k,w)\in
  \Pi_{G}^{\rm min}(v,w|C)$ then no non-consecutive vertices on $\pi$ are
  adjacent.
\end{lemma}
\begin{proof}
  Let $v_0=v$ and $v_{k+1}=w$, and suppose for a contradiction that there
  are non-consecutive vertices on the path $\pi$ which are adjacent. Let
  $(v_p,v_q)$ be a pair of adjacent vertices which are furthest apart on
  the path, i.e., $(p,q)$ maximizes the distance $|r-s|$ among pairs of
  indices of adjacent vertices $v_r$ and $v_s$ on the path.  Since $\pi$ is
  of minimal length, $v\neq v_p$ or $w\neq v_q$.

Suppose that $v\neq v_p$. By definition of $(p,q)$, $v_{p-1}$ is not
adjacent to $v_q$. Consequently, by Lemma \ref{lem:connect1}, $v_p$ is a
collider on $(v_{p-1},v_p,v_q)$, and thus the edge between $v_{p-1}$ and
$v_p$ has an arrowhead at $v_p$. It then follows by Lemma \ref{lem:tsr1}
that either $v_p\in \an(C)$ or $(v_p,v_{p+1},\ldots,v_k,w)$ is a directed
path from $v_p$ to $w$. In the latter case $v_p \in \an(v_q)$, but there is
an arrowhead at $v_p$ on the edge between $v_p$ and $v_q$, which
contradicts that $G$ is ancestral. Hence 
$v_p\in \an(C)$.  If $v_q=w$ then the path $(v,v_1,\ldots,v_p,v_q=w)$ is
$m$-connecting given $C$ and shorter than $\pi$. Hence $v_q \neq w$. It
then follows by the same argument that $v_q$ is a collider on
$(v_p,v_q,v_{q+1})$ and in $\an (C)$. However, this also leads to a
contradiction since then the path
$(v,v_1,\ldots,v_p,v_q,v_{q+1},\ldots,v_k,w)$ is both $m$-connecting given
$C$ and shorter than $\pi$.

The case where $w\neq v_q$ may be argued symmetrically.
\end{proof}

\begin{corollary}\label{cor:connect5}
  If $G$ is an ancestral graph that satisfies the boundary
  containment property and $\pi=(v=v_0,v_1,\ldots,v_k,v_{k+1}=w)\in
  \Pi_{G}^{\rm min}(v,w|C)$, then all the non-endpoint vertices
  $v_1,\ldots,v_k$ are colliders on $\pi$.
\end{corollary}
\begin{proof}
  This follows directly from Lemma \ref{lem:tsr2} and Lemma
  \ref{lem:connect1}.
\end{proof}

Even though all non-endpoints on a path of the type described in Corollary
\ref{cor:connect5} in $\Pi_{G}^{\rm min}(v,w|C)$ are colliders, not all
non-endpoints must be in the set $C$.  For example, in the graph $G^{\min}_2$
from Figure \ref{fig:Gs}, the path $(x,v,y)$ $m$-connects $x$ and $y$
given $\{ w\}$ since the collider $v$ is an ancestor of $w$. However, as
the next Lemma shows, there will always exist a path in $\Pi_{G}^{\rm
  min}(v,w|C)$ such that all non-endpoints are colliders in $C$.  In
$G^{\min}_2$ from Figure \ref{fig:Gs}, the path $(x,w,y)$ $m$-connects $x$
and $y$ given $\{w\}$.
\begin{lemma}\label{prop:connect6}
  If $G$ is an ancestral graph that satisfies the boundary
  containment property, and $\pi =
  (v=v_0,v_1,\ldots,v_k,v_{k+1}=w)\in\Pi_{G}^{\rm min}(v,w|C)$ is such that
  no other path in $\Pi_{G}^{\rm min}(v,w|C)$ has more non-endpoint
  vertices in $C$ than $\pi$, then all non-endpoint vertices
  $v_1,\ldots,v_k$ on $\pi$ are colliders that are in $C$.
\end{lemma}
\begin{proof} 
  By Corollary \ref{cor:connect5}, all non-endpoints $v_1,\ldots,v_k$ are
  colliders.  Assume that there exists $v_i\notin C$, $1\le i\le k$.  Since
  $\pi\in\Pi_{G}(v,w|C)$, and thus $v_i\in\an(C)$, there exists $c\in C$
  such that $v_i\to \cdots\to c\in G$. In particular, $c\not= v_{i-1}$ and
  $c\not= v_{i+1}$ because $v_i$ is ancestral neither to $v_{i-1}$ nor to
  $v_{i+1}$.  The boundary containment property and the fact that
  $G$ does not contain directed cycles imply that $v_i\to c\in G$.  By
  Lemma \ref{lem:connect1}, $G$ contains edges between $c$ and both
  $v_{i-1}$ and $v_{i+1}$.  Since the edge between $v_{i-1}$ and $v_i$ has
  an arrowhead at $v_i$ and $v_i\to c\in G$, the edge between $v_{i-1}$ and
  $c$ must have an arrowhead at $c$ because otherwise the fact that $G$ is
  an ancestral graph would be contradicted.  Similarly, the edge between
  $v_{i+1}$ and $c$ must have an arrowhead at $c$.  If $v_{i-1}\to c$, then
  $v_{i-2}$ is adjacent to $c$ and by the same argument as above there must
  be an arrowhead at $c$ on the edge between $v_{i-2}$ and $c$.  Repeating
  this argument yields that there exists a vertex $v_\ell$, $\ell\le i-1$,
  such that either $v_\ell\bi c\in G$, or $v_\ell=v$ and $v\to c$.  The
  same arguments also imply that there exists a vertex $v_j$, $j\ge i+1$,
  such that either $v_j\bi c\in G$, or $v_j=w$ and $w\to c$.  Therefore,
  the path $(v,v_1,\ldots,v_\ell,c,v_j,\ldots, v_k,w)$ is in
  $\Pi_{G}(v,w|C)$ and is either shorter than $\pi$ or of equal length but
  with more non-endpoint vertices in $C$.  This contradicts the choice of
  $\pi$ and therefore the assumption of a non-endpoint on $\pi$ that is not
  in $C$ must be false.
\end{proof}

\end{appendix}

\bibliographystyle{abbrvnat}
\bibliography{simplicial}

\end{document}